\input{epsf} 
\newcommand{\et}{\hspace{-0.08in}{\bf .}\hspace{0.1in}}
\newcommand{\BOX}{\hbox {$\sqcap$ \kern -1em $\sqcup$}}
\newcommand{\qed}{\hskip 3em \hbox{\BOX} \vskip 2ex}
\renewcommand{\hom}{{\rm hom}}
\newcommand{\bihom}{{\rm bihom}}

\newcommand{\To}{\Rightarrow}
\newcommand{\Set}{{\rm Set}}
\newcommand{\Cat}{{\rm Cat}}
\newcommand{\Rep}{{\rm Rep}}
\newcommand{\Cob}{{\rm Cob}}
\renewcommand{\to}{\rightarrow}
\newcommand{\tensor}{\otimes}
\newcommand{\maps}{\colon}
\newcommand{\op}{{\rm op}}
\newcommand{\eps}{\epsilon}
\newcommand{\iso}{\cong}
\newcommand{\cok}{{\rm cok}}
\newcommand{\im}{{\rm im}}

\newcommand{\Hilb}{{\rm Hilb}}
\newcommand{\SuperHilb}{{\rm SuperHilb}}
\renewcommand{\H}{H}
\newcommand{\Ab}{{\rm Ab}}
\newcommand{\Spec}{{\rm Spec}}
\newcommand{\tr}{{\rm tr}}
\newcommand{\qtr}{{\rm qtr}}
\newcommand{\qdim}{{\rm qdim}}

\newcommand{\SU}{{\rm SU}}
\newcommand{\U}{{\rm U}}
\renewcommand{\O}{{\rm O}}
\newcommand{\SO}{{\rm SO}}
\newcommand{\Sp}{{\rm Sp}}

\newcommand{\End}{{\rm end}}

\newcommand{\C}{{\Bbb C}}
\newcommand{\N}{{\Bbb N}}
\newcommand{\Z}{{\Bbb Z}}
\newtheorem{thm}{Theorem}    %[section]
\newtheorem{cor}[thm]{Corollary}

\newtheorem{prop}[thm]{Proposition}
\newtheorem{defn}[thm]{Definition}

        \newcommand{\be}{\begin{equation}}
        \newcommand{\ee}{\end{equation}}
        \newcommand{\ba}{\begin{eqnarray}}
        \newcommand{\ea}{\end{eqnarray}}
        \newcommand{\ban}{\begin{eqnarray*}}
        \newcommand{\ean}{\end{eqnarray*}}
        \newcommand{\barr}{\begin{array}} 
        \newcommand{\earr}{\end{array}}

\documentstyle[12pt,diagram,amsfonts]{article}
\begin{document}

      \begin{center}
      {\bf Higher-Dimensional Algebra II: \\
       2-Hilbert Spaces\\}
      \vspace{0.5cm}
      {\em John C.\ Baez\\}
      \vspace{0.3cm}
      {\small Department of Mathematics,  University of California\\ 
      Riverside, California 92521 \\
      USA\\ }
      \vspace{0.3cm}
      {\small email: baez@math.ucr.edu\\}
      \vspace{0.3cm}
      {\small August 27, 1996 \\ }
      \end{center}

\begin{abstract}

A 2-Hilbert space is a category with structures and properties
analogous to those of a Hilbert space.  More precisely, we define a
2-Hilbert space to be an abelian category enriched over $\Hilb$ with a
$\ast$-structure, conjugate-linear on the $\hom$-sets, satisfying
$\langle fg,h \rangle = \langle g,f^\ast h \rangle = \langle
f,hg^\ast\rangle$.  We also define monoidal, braided monoidal, and
symmetric monoidal versions of 2-Hilbert spaces, which we call
2-H*-algebras, braided 2-H*-algebras, and symmetric 2-H*-algebras, and
we describe the relation between these and tangles in 2, 3, and 4
dimensions, respectively.  We prove a generalized Doplicher-Roberts
theorem stating that every symmetric 2-H*-algebra is equivalent to the
category $\Rep(G)$ of continuous unitary finite-dimensional
representations of some compact supergroupoid $G$.  The equivalence is
given by a categorified version of
the Gelfand transform; we also construct a categorified version of the
Fourier transform when $G$ is a compact abelian group.  Finally, we
characterize $\Rep(G)$ by its universal properties when $G$ is a compact
classical group.  For example, $\Rep(\U(n))$ is the free connected symmetric
2-H*-algebra on one even object of dimension $n$.
\end{abstract}

\section{Introduction}

A common theme in higher-dimensional algebra is `categorification':
the formation of $(n+1)$-categorical analogs of $n$-categorical
algebraic structures.  This amounts to replacing equations between
$n$-morphisms by specified $(n+1)$-isomorphisms, in accord with the
philosophy that any {\it interesting} equation --- as opposed to one
of the form $x = x$ --- is better understood as an isomorphism, or
more generally an equivalence.

In their work on categorification in topological quantum field
theory, Freed \cite{Freed} and Crane \cite{Crane} have, in an
informal way, used the concept of a `2-Hilbert space': a category
with structures and properties analogous to those of a
Hilbert space.  Our goal here is to define 2-Hilbert spaces
precisely and begin to study them.  We concentrate on the
finite-dimensional case, as the infinite-dimensional case
introduces extra issues that we are not yet ready to handle.  We
must start by categorifying the various ingredients in the
definition of Hilbert space.  These are: 1) the zero element, 2)
addition, 3) subtraction, 4) scalar multiplication, and 5) the
inner product.  The first four have well-known categorical
analogs.

1) The analog of the zero vector is a `zero object'.  A zero
object in a category is an object that is both initial and
terminal.  That is, there is exactly one morphism from it to any
object, and exactly one morphism to it from any object.  Consider
for example the category $\Hilb$ having finite-dimensional
Hilbert spaces as objects, and linear maps between them as
morphisms.  In $\Hilb$, any zero-dimensional Hilbert space is a
zero object.

2) The analog of adding two vectors is forming the direct sum, or
more precisely the `coproduct', of two objects.  A coproduct of the
objects $x$ and $y$ is an object $x \oplus y$, equipped with
morphisms from $x$ and $y$ to it, that is universal with respect
to this property.  In $\Hilb$, for example, any Hilbert space equipped 
with an isomorphism to the direct sum of $x$ and $y$ is a coproduct of
$x$ and $y$.

3) The analog of subtracting vectors is forming the `cokernel' of
a morphism $f \maps x \to y$.  This makes sense only in a
category with a zero object.  A cokernel of $f \maps x \to y$
is an object $\cok f$ equipped with an epimorphism $g \maps y \to
\cok f$ for which the composite of $f$ and $g$ factors through
the zero object, that is universal with respect to this property. 
Note that while we can simply subtract a number $x$ from a number
$y$, to form a cokernel we need to say how the object $x$ is
mapped to the object $y$.  In $\Hilb$, for example, any space
equipped with an isomorphism to the orthogonal complement of $\im
f$ in $y$ is a cokernel of $f \maps x \to y$.  If $f$ is an
inclusion, so that $x$ is a subspace of $y$,  its cokernel is
sometimes written as the `direct difference' $y \ominus x$ to
emphasize the analogy with subtraction.

An important difference between zero, addition and subtraction and
their categorical analogs is that these operations represent extra
{\it structure} on a set, while having a zero object, binary
coproducts or cokernels is merely a {\it property} of a category.
Thus these concepts are in some sense more intrinsic to categories
than to sets.  On the other hand, one pays a price for this:
while the zero element, sums, and differences are unique in a Hilbert space,
the zero object, coproducts, and cokernels are typically 
unique only up to canonical isomorphism.   

4) The analog of multiplying a vector by a complex number is tensoring
an object by a Hilbert space.  Besides its
additive properties (zero object, binary coproducts, and cokernels),
$\Hilb$ also has a compatible multiplicative structure, that is,
tensor products and a unit object for the tensor product.  In other
words, $\Hilb$ is a `ring category', as defined by Laplaza and Kelly
\cite{Kelly2,Laplaza}.    We expect it to play a role in 2-Hilbert space
theory analogous to the role played by the ring $\C$ of complex numbers
in Hilbert space theory.  Thus we expect 2-Hilbert spaces to be `module
categories' over $\Hilb$, as defined by Kapranov and Voevodsky \cite{KV}.

An important part of our philosophy here is that $\C$ is the primordial
Hilbert space: the simplest one, upon which the rest are modelled.  By
analogy, we expect $\Hilb$ to be the primordial 2-Hilbert space.  This
is part of a general pattern pervading higher-dimensional algebra; for
example, there is a sense in which $n\Cat$ is the primordial
$(n+1)$-category. The real significance of this pattern remains somewhat
mysterious.

5) Finally, what is the categorification of the inner
product in a Hilbert space?  It appears to be the `$\hom$ functor'.  The
inner product in a Hilbert space $x$ is a bilinear map 
\[ \langle\,\cdot\,,\, \cdot\, \rangle \maps \overline x \times x \to \C
\] 
taking each pair of elements $v,w \in x$ to the inner product
$\langle v,w\rangle$.   Here $\overline x$ denotes the conjugate
of the Hilbert space  $x$.  Similarly, the $\hom$ functor in a
category $C$ is a bifunctor 
\[ \hom(\,\cdot\, , \,\cdot\,)\maps C^\op \times C \to \Set \] 
taking each pair of objects $c,d \in C$
to the set $\hom(c,d)$ of morphisms from $c$ to $d$.
This analogy clarifies the relation between 
category theory and quantum theory that is so important in
topological quantum field theory.  In quantum theory the inner product
$\langle v,w \rangle$ is a {\it number} representing the 
amplitude to pass from $v$ to $w$, while in category
theory $\hom(c,d)$ is a {\it set} of morphisms passing from $c$ to $d$.  

To understand this analogy better, note that 
any morphism $f \maps x \to y$ in $\Hilb$
can be turned around or `dualized' to obtain a morphism $f^\ast \maps y
\to x$.  The morphism $f^\ast$ is called the adjoint of $f$, and satisfies 
\[ \langle fv,w \rangle = \langle v,f^\ast w \rangle \] 
for all $v \in x$, $w \in y$.  
The ability to dualize morphisms in this way is crucial to
quantum theory.  For example, observables are represented by
self-adjoint morphisms, while symmetries are represented by unitary
morphisms, whose adjoint equals their inverse.  

The ability to dualize morphisms in $\Hilb$ makes this category
very different from the category $\Set$, in which the only
morphisms $f \maps x \to y$ admitting any natural sort of `dual'
are the invertible ones. There is, however, 
duals for certain noninvertible morphisms in $\Cat$ --- namely,
adjoint functors.  The functor $F^\ast \maps D \to C$ is
said to be a right adjoint of the functor  $F \maps C \to D$ if
there is a natural isomorphism
\[ \hom(Fc,d) \iso \hom(c,F^\ast d) \] 
for all $c \in C$, $d \in D$.  The analogy to adjoints of
operators between Hilbert spaces is clear.  Our main point here
is that that this analogy relies on the more fundamental analogy
between the inner product and the $\hom$ functor.  

One twist in the analogy between the inner product and the $\hom$
functor is that the inner product for a Hilbert space takes values in
$\C$.  Since we are treating $\Hilb$ as the categorification of $\C$, 
the $\hom$-functor for a 2-Hilbert space should take values in $\Hilb$
rather than $\Set$.  In technical terms \cite{Kelly}, this
suggests that a 2-Hilbert space should be enriched over $\Hilb$.

To summarize, we expect that a 2-Hilbert space should be
some sort of category with 1) a zero object, 2) binary coproducts, and 3)
cokernels, which is 4) a $\Hilb$-module and 5) enriched over $\Hilb$.  
However, we also need a categorical analog for the equation
\[            \langle v,w\rangle = \overline{\langle w,v \rangle} \]
satisfied by the inner product in a Hilbert space.  That is, for any two
objects $x,y$ in a 2-Hilbert space there should be a natural isomorphism 
\[          \hom(x,y) \iso \overline{\hom(y,x)}  \]
where $\overline{\hom(y,x)}$ is the complex conjugate of the Hilbert space
$\hom(y,x)$.  (The fact that objects in $\Hilb$ have complex conjugates
is a categorification of the fact that elements of $\C$ have complex
conjugates.)  This natural isomorphism should also satisfy some coherence
laws, which we describe in Section 2.  We put these ingredients together
and give a precise definition of 2-Hilbert spaces in Section 3.

Why bother categorifying the notion of Hilbert space?  As already
noted, one motivation comes from the study of topological quantum
field theories, or TQFTs.  In the introduction to this series
of papers \cite{BD}, we proposed that
$n$-dimensional unitary extended TQFTs should be treated as
$n$-functors from a certain $n$-category $n\Cob$ to a certain
$n$-category $n\Hilb$.  Roughly speaking, the $n$-category $n\Cob$
should have 0-dimensional manifolds as objects, 1-dimensional
cobordisms between these as morphisms, 2-dimensional cobordisms
between these as 2-morphisms, and so on up to dimension $n$.  The
$n$-category $n\Hilb$, on the other hand, should have `$n$-Hilbert
spaces' as objects, these being $(n-1)$-categories with structures
and properties analogous to those of Hilbert
spaces.  (Note that an ordinary Hilbert space is a `1-Hilbert space', and
is a 0-category, or set, with extra structures and properties.)

An eventual goal of this series is to develop the framework needed to
make these ideas precise.  This will require work both on $n$-categories in
general --- especially `weak' $n$-categories, which are poorly
understood for $n > 3$ --- and also on the particular $n$-categories
$n\Cob$ and $n\Hilb$.   One of the guiding lights of weak $n$-category
theory is the chart shown in Figure 1.  This describes `$k$-tuply
monoidal $n$-categories' --- that is, $(n+k)$-categories with only one
$j$-morphism for $j < k$.  The entries only correspond to theorems for
$n+k \le 3$, but there is evidence that the pattern continues for
arbitrarily large values of $n,k$.  Note in particular how as we descend
each column, the $n$-categories first acquire a `monoidal' or tensor
product structure, which then becomes increasingly `commutative' in
character with increasing $k$, stabilizing at $k = n+2$. 

\vskip 0.5em
\begin{center}
{\small
\begin{tabular}{|c|c|c|c|}  \hline
      & $n = 0$   & $n = 1$    & $n = 2$          \\     \hline
$k = 0$  & sets      & categories & 2-categories     \\     \hline
$k = 1$  & monoids   & monoidal   & monoidal         \\    
      &           & categories & 2-categories     \\     \hline
$k = 2$  &commutative& braided    & braided          \\      
      &  monoids  & monoidal   & monoidal         \\
      &           & categories & 2-categories     \\     \hline
$k = 3$  &`       & symmetric  & weakly involutory \\   
      &           & monoidal   & monoidal         \\
      &           & categories & 2-categories     \\     \hline
$k = 4$  &`'      & `'         &strongly involutory\\      
      &           &            & monoidal         \\
      &           &            & 2-categories     \\     \hline
$k = 5$  &`'      &`'          & `'               \\
      &           &            &                  \\ 
      &           &            &                  \\      \hline
\end{tabular}} \vskip 1em
1. The category-theoretic hierarchy: expected results
\end{center}
\vskip 0.5em 

At least in the low-dimensional cases examined so far, the
$n$-categories of interest in topological quantum field theory have
simple algebraic descriptions.  For example, knot theorists are familiar
with the category of framed oriented 1-dimensional cobordisms embedded
in $[0,1]^3$.  We would call these `1-tangles in 3 dimensions'.  They
form not merely a category, but a braided monoidal category.  In fact,
they form the `free braided monoidal category with duals on one object',
the object corresponding to the positively oriented point.  More
generally, we expect that $n$-tangles in $n+k$ dimensions form the `free
$k$-tuply monoidal $n$-category with duals on one object', $C_{n,k}$.
By its freeness, we should be able to obtain a representation of
$C_{n,k}$ in any $k$-tuply monoidal $n$-category with duals by
specifying a particular object therein.

When the codimension $k$ enters the stable range $k \ge n + 2$ we hope
to obtain the `free stable $n$-category with duals on one object',
$C_{n,\infty}$.  A unitary extended TQFT should be a representation of
this in $n\Hilb$.  If as expected $n\Hilb$ is a stable $n$-category
with duals, to specify a unitary extended TQFT would then simply be to
specify a particular $n$-Hilbert space.  More generally, we expect an
entire hierarchy of $k$-tuply monoidal $n$-Hilbert spaces in analogy
to the category-theoretic hierarchy, as shown in Figure 2.  We also
hope that an object in a $k$-tuply monoidal $n$-Hilbert space $H$
will determine a representation of $C_{n,k}$ in $H$,
and thus an invariant of $n$-tangles in $(n+k)$ dimensions.

\vskip 0.5em
\begin{center}
{\small
\begin{tabular}{|c|c|c|c|}  \hline
      & $n = 1$   & $n = 2$    & $n = 3$          \\     \hline
$k = 0$  & Hilbert    & 2-Hilbert & 3-Hilbert  \\ 
         & spaces     & spaces    & spaces     \\ \hline
$k = 1$  & H*-algebras   & 2-H*-algebras   & 3-H*-algebras   \\   \hline
$k = 2$  &commutative& braided    & braided          \\      
      & H*-algebras   & 2-H*-algebras   & 3-H*-algebras       \\  \hline
$k = 3$  &`'         & symmetric  & weakly involutory \\   
      &           & 2-H*-algebras   & 3-H*-algebras        \\ \hline
$k = 4$  &`'         & `'         &strongly involutory\\      
      &           &            & 3-H*-algebras        \\    \hline
$k = 5$  &`'         &`'          & `'               \\
      &           &            &                  \\ 
      &           &            &                  \\      \hline
\end{tabular}} \vskip 1em
2. The quantum-theoretic hierarchy: expected results
\end{center}
\vskip 0.5em 

We are far from proving general results along these lines!  However,
in Section 4 we sketch the structure of $2\Hilb$ as a strongly
involutory 3-H*-algebra, and in Section 5 we define 2-H*-algebras,
braided 2-H*-algebra, and symmetric 2-H*-algebras, and describe their
relationships to 1-tangles in 2, 3, and 4 dimensions, respectively.

An exciting fact about the quantum-theoretic hierarchy is that it
automatically subsumes various branches of representation theory.
2-H*-algebras arise naturally as categories of unitary
representations of certain Hopf algebras, or more generally `Hopf
algebroids', which are to groupoids as Hopf algebras are to
groups \cite{Lu}.  Braided 2-H*-algebras arise in a similar way from
certain quasitriangular Hopf algebroids --- for example, quantum
groups --- while symmetric 2-H*-algebras arise from certain
triangular Hopf algebroids --- for example, groups.

In Section \ref{recon} of this paper we concentrate on the
symmetric case.  Generalizing the Doplicher-Roberts theorem
\cite{DR}, we prove that all symmetric 2-H*-algebras are
equivalent to categories of representations of `compact
supergroupoids'.  If a symmetric 2-H*-algebra is `purely
bosonic', it is equivalent to a category of representations of a
compact groupoid; if it is `connected', it is equivalent to a
category of representations of a compact supergroup.  In
particular, any connected even symmetric 2-H*-algebra is
equivalent to the category $\Rep(G)$ of continuous
unitary finite-dimensional representations
of a compact group $G$.  This is the original Doplicher-Roberts
theorem.

One can view our generalized Doplicher-Roberts theorem as a
categorified version of the Gelfand-Naimark theorem.  The
Gelfand-Naimark theorem applies to commutative C*-algebras, but
one can easily deduce a version for commutative H*-algebras.
Roughly speaking, this says that every commutative H*-algebra
$H$ is isomorphic to a commutative H*-algebra of {\it
functions} from some set $\Spec(H)$ to
$\C$.  Similarly, our theorem implies that every even symmetric
2-H*-algebra $H$ is equivalent to a symmetric 2-H*-algebra of
{\it functors} from some groupoid
$\Spec(H)$ to $\Hilb$.  The equivalence is given explicitly by a
categorified version of the Gelfand transform.  We also construct
a categorified version of the Fourier transform, applicable to
the representation theory of compact abelian groups.

These links between the quantum-theoretic hierarchy and
representation theory give new insight into the
representation theory of classical groups.  The designation of a
group as `classical' is more a matter of
tradition than of some conceptual definition, but in practice
what makes a group `classical' is that it has a nice right
universal property.  In other words, there is a simple
description of homomorphisms into it.  Using the fact that group
homomorphisms from $G$ to $H$ determine symmetric 2-H*-algebra
homomorphisms from $\Rep(H)$ to $\Rep(G)$, one can show that for
a classical group $H$ the symmetric 2-H*-algebra $\Rep(H)$ has
nice left universal property: there is a simple description of
homomorphisms out of it.

For example, the group $\U(n)$ has a distinguished
$n$-dimensional unitary representation $\rho$, its fundamental
representation on $\C^n$.  An $n$-dimensional unitary
representation of any group $G$ is essentially the same as
a homomorphism from $G$ to $\U(n)$.  Using this right universal
property of $\U(n)$, we
show in Section \ref{recon} that the category of unitary
representations of $\U(n)$ is the `free symmetric 2-H*-algebra on
one object of dimension $n$'.  This statement tersely encodes the
usual description of the representations of $\U(n)$ in terms of
Young diagrams.  We also give similar characterizations of the
categories of representations of other classical groups.  

In what follows, we denote the composition of 1-morphisms, the
horizontal composition of a 1-morphism and a 2-morphism (in
either order) and the horizontal composition of 2-morphisms is
denoted by $\circ$ or simply juxtaposition.  Vertical composition
of 2-morphisms is denoted by $\cdot\,$.  {\it Nota bene}: in
composition we use the ordering in which, for example, the
composite of $f \maps x \to y$ and $g \maps y \to z$ is denoted
$f\circ g$.  We denote the identity morphism of an object $x$
either as $1_x$ or, if there is no danger of confusion, simply as
$x$.  We refer to our earlier papers on higher-dimensional
algebra as HDA0 \cite{BD} and HDA1 \cite{BN}.

\section{H*-Categories}

Let $\Hilb$ denote the category whose objects are 
finite-dimensional Hilbert spaces, and whose morphisms are arbitrary
linear maps.  (Henceforth, all Hilbert spaces will taken as
finite-dimensional unless otherwise specified.)  The category $\Hilb$
is symmetric monoidal, with $\C$ as the unit object, the usual tensor
product of Hilbert spaces as the monoidal structure, and the maps
\[            S_{x,y}(v \tensor w) = w \tensor v  \]
as the symmetry, where $x,y \in \Hilb$, $v \in x$, and $w \in y$.
Using enriched category theory \cite{Kelly} we may thus define
the notion of a category enriched over $\Hilb$, or $\Hilb$-category.
Concretely, this amounts to the following:

\begin{defn}\et  A $\Hilb$-{\rm category} $H$ is a category such that
for any pair of objects $x,y \in H$ the set of morphisms $\hom(x,y)$ is
equipped with the structure of a Hilbert space, and for any
objects $x,y,z \in \H$ the composition map 
\[     \circ\; \maps  \hom(x,y) \times \hom(y,z) \to \hom(x,z) \]
is bilinear.  \end{defn}

We may think of the `\hom' in a $\Hilb$-category $H$ as a functor
\[            \hom \maps \H^{\op} \times H \to \Hilb \]
as follows.  An object in $\H^\op \times H$ is just a pair of objects
$(x,y)$ in $\H$, and the $\hom$ functor assigns to this the object
$\hom(x,y) \in \Hilb$.  A morphism $F \maps (x,y) \to (x',y')$ in 
$\H^\op \times \H$ is just a pair of
morphisms $f \maps x' \to x$, $g \maps y \to y'$ in $\H$, and the
$\hom$ functor assigns to $F$ the morphism $\hom(F) \maps
\hom(x,y) \to \hom(x',y')$ given by
\[      \hom(F)(h) =  fhg .\]

As described in the introduction, we may regard $\Hilb$ as the
categorification of $\C$.  A structure on $\C$ which is crucial for
Hilbert space theory is complex conjugation,
\[        {}^{\overline{\hbox{\hskip 0.5em}}}\, \maps \C \to \C. \]
The categorification of this map is a functor
\[       {}^{\overline{\hbox{\hskip 0.5em}}}\, \maps \Hilb \to \Hilb \]
called {\it conjugation}, defined as follows.
First, for any Hilbert space $x$, there is a conjugate Hilbert space 
$\overline x$.  This has the same underlying abelian
group as $x$, but to keep things straight let us temporarily write
$\overline v$ for the element of $\overline x$
corresponding to $v \in x$.  Scalar multiplication in $\overline x$ is
then given by 
\[            c\overline v = \overline {(\overline c v)} \]
for any $c \in \C$, while the inner product is given by 
\[      \langle \overline v, \overline w \rangle = \overline {\langle
v,w\rangle}. \]
Second, for any morphism $f \maps x \to y$ in $\Hilb$, there is
a conjugate morphism $\overline f \maps \overline x \to \overline y$,
given by 
\[          \overline f(\overline v) = \overline{f(v)} \]
for all $v \in x$.  One can easily check that with these
definitions conjugation is a covariant functor. Note that the
square of this functor is equal to the identity.   Also note that
a linear map $f \maps x \to \overline y$ is the same thing as an
antilinear (i.e., conjugate-linear) map from $x$ to $y$, while a
unitary map $f \maps x \to \overline y$ is the same thing as an
antiunitary map from $x$ to $y$.

Now, just as in a Hilbert space we have the equation
\[        \langle v,w \rangle = \overline {\langle w,v\rangle} \]
for any pair of elements, in a 2-Hilbert space we would like
an isomorphism 
\[        \hom(x,y) \iso \overline {\hom(y,x) } \]
for every pair of objects.  This isomorphism should be be
`natural' in some sense, but $\hom(x,y)$ is contravariant in $x$
and covariant in $y$, while $\overline {\hom(y,x)}$ is covariant
in $x$ and contravariant in $y$.   Luckily $\Hilb$ is a
$\ast$-category, which allows us to define `antinatural
isomorphisms' between covariant functors and  contravariant
functors from any category to $\Hilb$.  

This works as follows.  In general, a {\it $\ast$-structure} for a
category $C$ is defined as a contravariant functor $\ast \maps C \to C$
which acts as the identity on the objects of $C$ and satisfies $\ast^2 =
1_C$.  A {\it $\ast$-category} is a category equipped with a
$\ast$-structure.   For example, $\Hilb$ is a $\ast$-category 
where for any morphism $f\maps x \to y$ we define $f^\ast \maps y \to x$
to be the Hilbert space adjoint of $f$: 
\[            \langle fv, w\rangle = \langle v, f^\ast w \rangle \]
for all $v \in x$, $w \in y$.  

Now suppose that $D$ is a $\ast$-category and $F \maps C \to D$ is a
covariant functor, while $G \maps C \to D$ is a contravariant functor.  
We define an {\it antinatural transformation} $\alpha \maps F \To G$
to be a natural transformation 
from $F$ to $G \circ \ast$.  
Similarly, an antinatural transformation from $G$ to $F$ is defined to 
be a natural transformation from $G$ to $F \circ \ast$.

As a step towards defining a 2-Hilbert space we now define an
H*-category.
\begin{defn}\et An {\rm H*-category} is a $\Hilb$-category with a
$\ast$-structure that defines an antinatural transformation from
$\hom(x,y)$ to $\overline {\hom(y,x)}$.
\end{defn}

\noindent
This may require some clarification.  
Given a $\ast$-structure $\ast
\maps \H \to 
\H$, we obtain for any objects $x,y \in \H$ a function $\ast \maps \hom(x,y)
\to \hom(y,x)$.  By abuse of notation we may also regard this as a function
\[       \ast \maps \hom(x,y) \to \overline{\hom(y,x)}. \]
We then demand that this define an antinatural
transformation between the covariant functor $\hom \maps H^\op \times
H \to \Hilb$ to the contravariant functor sending $(x,y) \in
H^\op \times H$ to $\overline{\hom(y,x)} \in \Hilb$.  

The following proposition gives a more concrete description of
H*-categories: 
\begin{prop}\label{H*1}\et An H*-category $\H$ is the same as a
$\Hilb$-category equipped with antilinear maps $\ast \maps \hom(x,y) \to
\hom(y,x)$ for all $x,y \in \H$, such that 
\begin{enumerate} 
\item $f^{\ast\ast} = f$, 
\item $(fg)^\ast = g^\ast f^\ast$, 
\item $\langle fg,h\rangle = \langle g,f^\ast h \rangle$, 
\item $\langle fg,h\rangle = \langle f, hg^\ast \rangle$
\end{enumerate}
whenever both sides of the equation are well-defined.
\end{prop}

Proof - First suppose that $\H$ is an H*-category.  By the
antinaturality of $\ast$, for all $x,y \in \H$ there is a linear
map $\ast \maps \hom(x,y) \to \overline{\hom(y,x)}$, which is the same as
an antilinear map $\ast \maps \hom(x,y) \to \hom(y,x)$.   The fact that
$\ast$ is a $\ast$-structure implies properties 1 and 2.  
As for 3 and 4, suppose $(x,y)$ and 
$(x',y')$ are objects in $\H^\op \times \H$, and let $(f,g)$ be a
morphism from $(x,y)$ to $(x',y')$.  
The fact that $\ast$ is an antinatural transformation means that
the following diagram commutes:
\[ 
\begin{diagram}[\hom(x'y')]
\node{\hom(x,y)} \arrow{e,t}{\ast} \arrow{s,l}{V(f,g)} 
\node{\overline{\hom(y,x)}} \arrow{s,r}{W(f,g)^\ast}       \\
\node{\hom(x',y')}  \arrow{e,t}{\ast}    \node{\overline{\hom(y',x')}}
\end{diagram}
\]
where $V$ is the covariant functor 
\[
\begin{diagram}[xxxxxxx]
\node{\H^\op \times \H}
\arrow{e,t}{\hom}
\node{\Hilb}
\end{diagram}
\]
and $W$ is the contravariant functor
\[
\begin{diagram}[xxxxxxx]
\node{\H^\op \times \H}
\arrow{e,t}{S_{\H^\op,\H}}
\node{\H \times \H^\op}
\arrow{e,t}{\hom}
\node{\Hilb}
\arrow{e,t}{{}^{\overline{\hbox{\hskip 0.5em}}}}
\node{\Hilb,}
\end{diagram}
\]
where in this latter diagram
$S$ denotes the symmetry in $\Cat$, $\hom$ is regarded as a
contravariant functor from $(\H^\op \times \H)^\op \iso
\H \times \H^\op$ to $\Hilb$, and the overline denotes conjugation.
This is true if and only if for all $h \in \hom(x,y)$ 
and $k \in \overline{\hom(y',x')}$, 
\[     \langle (V(f,g)h)^\ast ,k \rangle = 
\langle W(f,g)^\ast h^\ast, k\rangle \]
or in other words, 
\[       \langle (fhg)^\ast, k\rangle = \langle h^\ast, gkf\rangle \]
or
\[   \langle  g^\ast h^\ast f^\ast, k\rangle = \langle h^\ast, gkf\rangle
.\]
Here the inner products are taken in $\overline{\hom(y',x')}$, but the
equations also hold with the inner product taken in $\hom(y',x')$.  
Taking either $f$ or $g$ to be the identity, we obtain
3 and 4 after some relabelling of variables.   

Conversely, given antilinear maps $\ast \maps \hom(x,y) \to
\hom(y,x)$ for all $x,y \in \H$, properties 1 and 2 say that these
define a $\ast$-structure for $\H$, and using 3 and 4 we obtain
\ban   \langle g^\ast h^\ast f^\ast , k\rangle &=&
  \langle g^\ast h^\ast , kf\rangle \\
&=&  \langle  h^\ast, gkf\rangle ,\ean
showing that $\ast$ is antinatural.  
\qed 

\begin{cor}\label{H*2}\et If $H$ is an H*-category, for all objects $x,y
\in H$ the map $\ast \maps \hom(x,y) \to \hom(y,x)$ is antiunitary.
\end{cor}

Proof - The map $\ast \maps \hom(x,y) \to \hom(y,x)$ is antilinear, and
by 3 and 4 of Proposition \ref{H*1} we have
\[         \langle f,g \rangle = \langle g^\ast, f^\ast \rangle 
= \overline{\langle f^\ast, g^\ast \rangle} \]
for all $f,g \in \hom(x,y)$, so $\ast$ is antiunitary.  \qed

Next we give a structure theorem for H*-categories.  This relies
heavily on the theory of `H*-algebras' due to Ambrose
\cite{Ambrose}, so let us first recall this theory.   For our
convenience, we use a somewhat different definition of
H*-algebra than that given by Ambrose.  Namely, we restrict our
attention to finite-dimensional H*-algebras with multiplicative unit,
and we do not require the inequality $\|ab\| \le \|a\| \, \|b\|$.

\begin{defn} \et An {\it H*-algebra} $A$ is a Hilbert space
that is also an associative algebra with unit, equipped with
an antilinear involution $\ast \maps A \to A$ satisfying 
\ban     
\langle ab,c\rangle &=& \langle b,a^\ast c \rangle       \\
\langle ab,c\rangle &=& \langle a, cb^\ast \rangle   \ean
for all $a,b,c \in A$.
An {\rm isomorphism} of H*-algebras is a unitary operator
that is also an involution-preserving algebra isomorphism.  
\end{defn} 

The basic example of an H*-algebra is the space of 
linear operators on a Hilbert space $H$.  
Here the product is the usual product of operators, the involution is
the usual adjoint of operators, and the inner product is given by
\[         \langle a,b\rangle = k \,\tr(a^\ast b) \]
where $k > 0$.  We denote this H*-algebra by $L^2(H,k)$. 
It follows from the work of Ambrose that all H*-algebras can be built
out of H*-algebras of this form.  More precisely, every H*-algebra $A$
is the orthogonal direct sum of finitely many 
minimal 2-sided ideals $I_i$, each
of which is isomorphic as an H*-algebra to $L^2(H_i,k_i)$ 
for some Hilbert space $H_i$ and some positive real number $k_i$.  

This result immediately classifies H*-categories with one object.  Given
an H*-category with one object $x$, $\End(x)$ is an H*-algebra, and is
thus of the above form.  Conversely, any H*-algebra is isomorphic to
$\End(x)$ for some H*-category with one object $x$.

We generalize this to arbitrary H*-categories as follows.
Suppose first that $H$ is an H*-category with finitely many objects.
Let $A$ denote the orthogonal direct sum
\[            A = \bigoplus_{x,y} \hom(x,y)  .\]
Then $A$ becomes an H*-algebra if we define
the product in $A$ of morphisms in $H$ to be their composite when the
composite exists, and zero otherwise, and define the involution in $A$
using the $\ast$-structure of $H$.  $A$ is thus the
orthogonal direct sum of finitely many minimal 2-sided ideals:
\[         A = \bigoplus_{i = 1}^n L^2(H_i, k_i)
.\]
For each object $x \in H$, the identity morphism $1_x$ can be regarded as an
element of $A$.  This element is a self-adjoint projection, meaning that
\[       1_x^\ast = 1_x, \qquad 1_x^2 = 1_x .\]
It follows that we may write
\[          1_x = \bigoplus_{i = 1}^n p_i^x \]
where $p_i^x \in L^2(H_i, k_i)$ is the 
projection onto some subspace $H_i^x \subseteq H_i$.  
Note that the elements $1_x$, $x \in H$, form a complete orthogonal set of
projections in $A$.  In other words, $1_x 1_y  = 0$ if $x \ne y$, 
and
\[           \sum_{x \in H} 1_x = 1.\]
Thus each Hilbert space $H_i$ is the orthogonal direct sum of the
subspaces $H_i^x$.  

This gives the following structure theorem for H*-categories:
\begin{thm}\label{H*3}\et Let $H$ be an H*-category and $S$ any finite
set of objects of $H$.  Then for some $n$, there
exist positive numbers $k_i > 0$ and
Hilbert spaces $H_i^x$ for $i = 1, \dots, n$ and $x \in S$, 
such that the following hold:
\begin{enumerate}
\item {\rm For $i = 1, \dots, n$, let 
\[       H_i = \bigoplus_{x \in S} H_i^x \]
denote the orthogonal direct sum, and let $p_i^x$ be the
self-adjoint projection from $H_i$ to $H_i^x$.  Then
for any objects $x,y \in S$, there is a unitary isomorphism
between the Hilbert space $\hom(x,y)$ and the subspace
\[        \bigoplus_i p_i^x
L^2(H_i,k_i) p_i^y \; \subseteq \;  \bigoplus_i L^2(H_i, k_i) ,\] 
Thus we may write any morphism $f \maps x \to y$ as
\[         f =  \bigoplus_i f_i \]
where $f_i \maps H_i^x \to H_i^y$.  
\item Via the above isomorphism, the composition map
\[        \circ\;\maps \hom(x,y) \times \hom(y,z) \to \hom(x,z) \]
is given by
\[         f \circ g = \bigoplus_i f_i g_i. \]
\item Via the same isomorphism, the $\ast$-structure
\[         \ast \maps \hom(x,y) \to \hom(y,x)  \]
is given by
\[        f^\ast = \bigoplus_i f_i^\ast .\]
\rm}
\end{enumerate}
\noindent Conversely, given a $\Hilb$-category $H$ with $\ast$-structure
such that the above holds for any finite subset $S$ of its objects, $H$
is an H*-category.  
\end{thm}
Proof - If $H$ has finitely many objects and we take $S$ to be the set
of all objects of $H$, properties 1-3 follow from the remarks preceding
the theorem.  More generally, by Proposition \ref{H*1} any full subcategory
of an H*-category is an H*-category, so 1-3 hold for any finite subset $S$
of the objects of $H$.  

Conversely, given a $\Hilb$-category $H$ with a $\ast$-structure, if
every full subcategory of $H$ with finitely many objects is an
H*-category, then $H$ itself is an H*-category.  One may check using
Proposition \ref{H*1} that if $S$ is any finite subset of the objects
of $H$, properties 1-3 imply the full subcategory of $H$ with $S$ as
its set of objects is an H*-category.  Thus $H$ is an H*-category.
\qed

The notions of unitarity and self-adjointness
will be important in all that follows.  
\begin{defn} \et  Let $x$ and $y$ be objects of a $\ast$-category.  
A morphism $u \maps x \to y$ is {\rm unitary} if $uu^\ast = 1_x$ and
$u^\ast u = 1_y$.   A morphism $a \maps x \to x$ is {\rm
self-adjoint} if $a^\ast = a$.    \end{defn}
Note that every unitary morphism is an isomorphism.  Conversely, the
following proposition implies that in an H*-category, isomorphic objects
are isomorphic by a unitary.
\begin{prop} \label{H*4} \et Suppose $f \maps x \to y$ is an isomorphism in 
the H*-category $H$.  Then $f = au$ where $a \maps x \to x$ is self-adjoint
and $u \maps x \to y$ is unitary.
\end{prop} 
Proof - Suppose that $f \maps x \to y$ is an isomorphism.  Then applying
Theorem \ref{H*3} to the full subcategory of $H$ with $x$ and $y$ as its
only objects, we have $f = \bigoplus f_i$ with $f_i \maps H_i^x \to
H_i^y$ an isomorphism for all $i$.    Using the polar
decomposition theorem we may write $f_i = a_i u_i$, where $a_i \maps
H_i^x \to H_i^x$ is the positive square root of $f_i {f_i}^\ast$, and $u_i
\maps H_i^x \to H_i^y$ is a unitary operator given by $u_i = 
a_i^{-1} f_i$.     Then defining $a = \bigoplus a_i$ and $u = \bigoplus
u_i$, we have $f = au$ where $a$ is self-adjoint and $u$ is unitary. \qed

One can prove a more general polar decomposition theorem allowing one
to write any morphism $f \maps x \to y$ in an H*-category as the
product of a self-adjoint morphism $a \maps x \to x$ and a partial
isometry $i \maps x \to y$, that is, a morphism for which $ii^\ast$
and $i^\ast i$ are self-adjoint idempotents.  However, we will not
need this result here.

\section{2-Hilbert Spaces}

The notion of 2-Hilbert space is intended to be the categorification of
the notion of Hilbert space.  As such, it should be a category having a
zero object, direct sums and `direct differences' of objects, tensor
products of Hilbert spaces with objects, and `inner products' of objects.
So far, with our definition of H*-category, we have formalized the
notion of a category in which the `inner product' $\hom(x,y)$ of any
two objects $x$ and $y$ is a Hilbert space.  Now we deal with the rest
of the properties:

\begin{defn}\et A {\rm 2-Hilbert space} is an abelian H*-category.
\end{defn} 

Recall that an abelian category is an $\Ab$-category (a category
enriched over the category $\Ab$ of abelian groups) such that
\begin{enumerate}
\item There exists an initial and terminal object.
\item Any pair of objects has a biproduct.
\item Every morphism has a kernel and cokernel.
\item Every monomorphism is a kernel, and every epimorphism is a
cokernel. 
\end{enumerate}
Let us comment a bit on what this amounts to.  
Since an H*-category is enriched over $\Hilb$ it is automatically
enriched over $\Ab$.   We call an initial and terminal object a
{\it zero object}, 
and denote it by $0$.  The zero object in a 2-Hilbert space is the
analog of the zero vector in a Hilbert space.  We call the biproduct
of $x$ and $y$ the {\it direct sum},
and denote it by $x \oplus y$.  Recall that by definition, this is equipped 
with morphisms $p_x \maps x \oplus y \to x$, $p_x \maps x \oplus y \to y$, 
$i_x \maps x \to x \oplus y$, $i_y \maps y \to x \oplus y$ such 
that 
\[       i_x p_x = 1_x , \qquad i_y p_y = 1_y, \qquad
  p_x i_x + p_y i_y = 1_{x\oplus y}  .\]
The direct sums in a 2-Hilbert space are the analog of addition
in a Hilbert space.   Similarly, 
the cokernels in a 2-Hilbert space are the analogs of differences
in a Hilbert space.  Finally, the
ability to tensor objects in a 2-Hilbert space by Hilbert
spaces (the analog of scalar multiplication) will follow from the other
properties, so we do not need to include it in the definition of
2-Hilbert space.

Some aspects of our definition of 2-Hilbert space may seem unmotivated
by the analogy with Hilbert spaces.  Why should a 2-Hilbert space have
kernels, and why should it satisfy clause 4 in the definition of abelian
category?  In fact, these properties follow from the rest.  

\begin{prop}\et \label{2hilb1} Let $H$ be an H*-category.  Then
the following are equivalent:
\begin{enumerate}
\item There exists an initial object.
\item There exists a terminal object. 
\item There exists a zero object.
\end{enumerate}
Moreover, the following are equivalent:
\begin{enumerate}
\item Every pair of objects has a product. 
\item Every pair of objects has a coproduct.
\item Every pair of objects has a direct sum.
\end{enumerate}
Moreover, the following are equivalent: 
\begin{enumerate}
\item Every morphism has a kernel.
\item Every morphism has a cokernel.
\end{enumerate}
Finally, if $H$ has a zero object, every pair of objects in $H$ has
a direct sum, and every morphism in $H$ has a cokernel, then 
$H$ is a semisimple abelian category.  
\end{prop}
Proof - It is well-known \cite{MacLane} that an initial or terminal
object in an $\Ab$-category is automatically a zero object.
Alternatively, this is true in every $\ast$-category, using the
bijection $\ast \maps \hom(x,y) \to \hom(y,x)$.  
It is also well-known that in an $\Ab$-category,
a binary product or coproduct is automatically a binary biproduct.
Furthermore, it is easy to check that
in any $\ast$-category, the morphism $j \maps k \to x$ is a kernel of
$f\maps x \to y$ if and only if $j^\ast \maps x \to k$ is a cokernel of
$f^\ast \maps y \to x$.  Thus a $\ast$-category has kernels if and only
if it has cokernels.  

Now suppose that $H$ is an H*-category with a zero object, direct
sums, and cokernels.  Then $H$ has kernels as well, so to show $H$ is
abelian we merely need to prove that every monomorphism is a kernel and
every epimorphism is a cokernel.  Let us show a monomorphism $f \maps x
\to y$ is a kernel; it follows using the $\ast$-structure that every
epimorphism is a cokernel.  

It suffices to show this result for any full subcategory of
$H$ with finitely many objects, so by Theorem \ref{H*3}
we may write
\[     f = \bigoplus_i f_i \]
where $f_i \maps H_i^x \to H_i^y$ is a linear operator.    Let $p
\maps y \to y$ be given by $\bigoplus p_i$ where $p_i$ is the
projection onto the orthogonal complement of the range of $f_i$. 
We claim that $f\maps x \to y$ is a kernel of $p$.   Since $f_i
p_i = 0$ for all $i$ we have $fp  = 0$.  We also need to show
that if $f' \maps x' \to y$ is any morphism with $f'p = 0$, then
there is a unique $g \maps x' \to x$ with $f' = gf$. Writing $f'
= \bigoplus f'_i$, the fact that $f'p = 0$ implies that the range
of $f'_i$ is contained in the range of $f_i$.  Thus by linear
algebra there exists $g_i \maps H_i^{x'} \to H_i^x$ such that 
$f'_i = g_i f_i$.  Letting $g = \bigoplus g_i$, we have $f' =
gf$, and $g$ is  unique with this property because $f$ is monic.

Finally, note that $H$ is semisimple, i.e., every
short exact sequence splits.  This follows from Theorem \ref{H*3}
and elementary linear algebra.  \qed

Given a 2-Hilbert space $H$, the fact that $H$ is semisimple implies
that every object is isomorphic to a direct sum of {\it simple} objects,
that is, objects $x$ for which $\End(x)$ is isomorphic as an algebra
to $\C$.  This fact lets us reason about 2-Hilbert spaces using bases:
\begin{defn} \et  Given a 2-Hilbert space $H$, a set of
nonisomorphic simple objects of $H$ is called a {\rm basis} if every
object of $H$ is isomorphic to a finite direct sum of objects 
in that set.  \end{defn}

\begin{cor} \label{2hilb2} \et  
Every 2-Hilbert space $H$ has a basis, and any two bases of $H$
have the same cardinality.  \end{cor}
Proof - The 2-Hilbert space $H$ has a basis because it is
semisimple: given any  Given two bases $\{e_\alpha\}$ and $\{f_\beta\}$,
each object $e_\alpha$ is isomorphic to a direct sum of copies of
the objects $e_\beta$, but as the $e_\alpha$ and $f_\beta$ are
simple we must actually have an isomorphism $e_\alpha \iso
f_\beta$ for some $\beta$. This $\beta$ is unique since no
distinct $f_\beta$'s are isomorphic.  This sets up a function
from $\{e_\alpha\}$ to $\{f_\beta\}$, and similar reasoning gives
us the inverse function.  \qed
\begin{defn} \et  The {\rm dimension} of a 2-Hilbert space is the
cardinality of any basis of it.  \end{defn}

Note that every basis $\{e_\alpha\}$ of a 2-Hilbert space
is `orthogonal' in the sense that
\[ \hom(e_\alpha, e_\beta) \iso 
\cases{L^2(\C,k_\alpha) & $\alpha = \beta$ \cr
0 & $\alpha \ne \beta$ } \]
where the isomorphism is one of H*-algebras, 
and $k_\alpha$ are certain positive constants.  
Moreover, up to reordering, the constants $k_\alpha$ are independent of
the choice of basis.  For 
suppose $x,y$ are two isomorphic objects in an H*-category.  By 
Proposition \ref{H*4} there is a unitary isomorphism $f \maps x \to y$.
Then there is an H*-algebra isomorphism $\alpha \maps \End(x) \to
\End(y)$ given by $\alpha(g) = f^{-1}gf$.  

One would also like to be able to tensor objects in a 2-Hilbert space
with Hilbert spaces, but this is a consequence of the definition we
have given, since one may define the tensor product of an object $x$
in a 2-Hilbert space with an $n$-dimensional Hilbert space to be the
direct sum of $n$ copies of $x$.  In fact, $\Hilb$ has a structure
analogous to that of an algebra, with tensor product and direct sum
playing the roles of multiplication and addition.  In the terminology
we introduce in Section \ref{2-H*-algebras}, one says that $\Hilb$ is
a `2-H*-algebra'.  One can develop a theory of modules of
2-H*-algebras following the ideas of Kapranov and Voevodsky \cite{KV}
and Yetter \cite{Yetter}.  Every 2-Hilbert space $H$ is then a module
over $\Hilb$.  We will not pursue this further here.

\section{2Hilb as a 2-Category}\label{2Hilb}

We now investigate a certain 2-category $2\Hilb$ of 2-Hilbert spaces.  
To keep things simple we take as its objects only finite-dimensional
2-Hilbert spaces.  Nonetheless we prove theorems more generally whenever
possible.  

\begin{defn} \et A {\rm morphism} $F \maps H \to H'$ between 2-Hilbert
spaces $H$ and $H'$ is an exact functor such that
$F \maps \hom(x,y) \to \hom(F(x), F(y))$ is linear and $F(f^\ast) =
F(f)^\ast$ for all $f \in \hom(x,y)$. \end{defn}
Recall that an exact functor is one preserving short exact
sequences.  Exactness is an natural sort of condition for functors 
between abelian categories.  Similarly, the requirement that 
$F\maps \hom(x,y) \to \hom(F(x), F(y))$ be linear is a natural condition 
for functors between
$\Hilb$-categories; one calls such a functor a {\it $\Hilb$-functor}.
Finally, $F(f^\ast) = F(f)^\ast$ is a natural condition
for functors between $\ast$-categories, and functors satisfying it are 
called {\it $\ast$-functors}.  

The following fact is occaisionally handy:
\begin{prop} \et Let $F \maps H \to H'$ be a functor between 2-Hilbert
spaces such that for all $x,y \in H$, 
$F \maps \hom(x,y) \to \hom(F(x), F(y))$ is linear.
Then the following are equivalent:
\begin{enumerate}
\item $F$ is exact.
\item $F$ is left exact.
\item $F$ is right exact.
\item $F$ preserves direct sums.  
\end{enumerate} 
\end{prop}
Proof - Following Yetter \cite{Yetter}, we use the
fact that every short exact sequence splits. \hbox{\hskip 30em} \qed

\begin{defn} \et A {\rm 2-morphism} $\alpha \maps F \To F'$ between 
morphisms $F,F' \maps H \to H'$ between 2-Hilbert spaces $H$ and $H'$ is
a natural transformation.  \end{defn}
\begin{defn} \et
We define the 2-category $2\Hilb$ to be that for which
 objects are finite-dimensional 2-Hilbert spaces, while morphisms
and 2-morphisms are defined as above.  
\end{defn}

Now, just as in some sense $\C$ is the primordial Hilbert space and
$\Hilb$ is the primordial 2-Hilbert space, $2\Hilb$ should be the
primordial 3-Hilbert space.  The study of $2\Hilb$ should thus shed
light on the properties of the still poorly understood 3-Hilbert spaces.
However, note that $\C$ is not merely a Hilbert space, but also a
commutative monoid, in fact a commutative H*-algebra.  Similarly,
$\Hilb$ is not merely a 2-Hilbert space, but also a symmetric monoidal
category when equipped with its usual tensor product.  Indeed, in
Section \ref{2-H*-algebras} we show that $\Hilb$ is a `symmetric
2-H*-algebra'.  Likewise, we expect $2\Hilb$ to be not only a
3-Hilbert space, but also a strongly involutory monoidal 2-category, in
fact a `strongly involutory 3-H*-algebra'.  As sketched in HDA0,
commutative monoids, symmetric monoidal categories, and strongly
involutory monoidal 2-categories are all examples of `stable'
$n$-categories.  In general we expect $n\Hilb$ to be a `stable
$(n+1)$-H*-algebra.'  The results below offer some support for this
expectation.

We begin with a study of duality in $2\Hilb$, as this is the most
distinctive aspect of Hilbert space theory.  Note that every element $x
\in \C$ has a kind of `dual' element, namely, its complex conjugate
$\overline x$.  Similarly, the category $\Hilb$ has duality both for
objects and for morphisms.  At the level of morphisms, each linear map
$f \maps x \to y$ between Hilbert spaces has a dual $f^\ast \maps y \to
x$, the usual Hilbert space adjoint of $f$.  This defines a
$\ast$-structure on $H$.  Duality at the level of objects can be
regarded either as a contravariant functor assigning to each each
Hilbert space $x$ its dual $x^\ast$, or as a covariant functor assigning
to each Hilbert space $x$ its conjugate $\overline x$.  These
two viewpoints become equivalent if we take advantage of duality at the
morphism level, since $x^\ast$ and $\overline x$ are antinaturally
isomorphic. 

Similarly, $2\Hilb$ has duality for objects, morphisms, and 2-morphisms.
As in $\Hilb$, we can use duality at a given level to reinterpret
dualities at lower levels in various ways.  This recursive process can
become rather confusing unless we choose by convention to take certain
dualities as `basic' and others as derived.  Here we follow the
philosophy of HDA0: any 2-morphism $\alpha \maps F \To G$ has a dual
$\alpha^\ast \maps G \To F$, any morphism $F \maps H \to H'$ has a dual
$F^\ast \maps H' \to H$, and every object $H$ has a dual $H^\ast$.
(Our notation differs from HDA0 in that we use the same symbol
to denote all these different levels of duality.)

\subsection{Duality for 2-morphisms}

Duals of 2-morphisms are the easiest to define.  It pays to do so in
the greatest possible generality:
\begin{defn} \label{star} \et Given a category $C$ and a $\ast$-category $D$,
the {\rm dual} $\alpha^\ast$ of a natural transformation $\alpha \maps F
\To G$ is the natural transformation with $(\alpha^\ast)_c =
(\alpha_c)^\ast$ for all $c \in C$.  \end{defn}
It is easy to check that $\alpha^\ast$ is a natural transformation when
$\alpha$ is, and that
\[        (\alpha^\ast)^\ast = \alpha, \qquad 1^\ast = 1.\]
The vertical composite of natural transformations satisfies 
\[     (\alpha\cdot \beta)^\ast = \beta^\ast\cdot \alpha^\ast \]
when this is defined.  When $D$ is a $\ast$-category, the horizontal
composite of a functor $F \maps B \to C$ and a natural
transformation $\alpha \maps G \To H$ with $G,H \maps C \to D$ satisfies
\[       (F\alpha)^\ast = F \alpha^\ast .\]
Similarly, when $F \maps C \to D$ is a $\ast$-functor and $\alpha
\maps G \To H$ is a natural transformation between $G,H \maps B
\to C$, we have
\[        (\alpha F)^\ast = \alpha^\ast F.\]

In particular, taking $C,D$ to be 2-Hilbert spaces, we obtain 
the definition of the dual of a 2-morphism in $2\Hilb$.  
We also obtain the notion of `unitary' and `self-adjoint' natural
transformations: 
\begin{defn} \et 
Given a category $C$, a $\ast$-category $D$, and functors $F,G \maps C
\to D$, a natural transformation $\alpha \maps F \To G$ is {\rm unitary} if
\[   \alpha \alpha^\ast = 1_F , \qquad \alpha^\ast \alpha = 1_G .\]
A natural transformation $\alpha \maps F \To F$ is {\rm self-adjoint} if
\[       \alpha^\ast = \alpha .\]
\end{defn}
Equivalently, $\alpha$ is unitary if $\alpha_c$ is a unitary morphism in
$D$ for all objects $c \in C$, and self-adjoint if $\alpha_c$ is
self-adjoint for all $c \in C$.  

Note that every unitary natural transformation is a natural isomorphism.
Conversely: 
\begin{prop} \label{unitary.nat} \et Suppose $F,G \maps H \to H'$ are
morphisms between 2-Hilbert spaces and $\alpha \maps F \To G$ is a natural
isomorphism.   Then $\alpha = \beta\cdot\gamma$ where $\beta \maps F \To
F$ is self-adjoint and $\gamma \maps F \To G$ is unitary.  \end{prop}

Proof - By Proposition \ref{H*4}, for any $x \in H$ we can write the
isomorphism $\alpha_x \maps F(x) \to G(x)$ as the composite $\beta_x
\gamma_x$, where $\beta_x \maps F(x) \to F(x)$ is self-adjoint and $\gamma_x
\maps F(x) \to G(x)$ is unitary.  More importantly, the 
polar decomposition gives a
natural way to construct $\beta_x$ and $\gamma_x$ from $\alpha_x$: we
take $\beta_x$ to be the positive square root of $\alpha_x
{\alpha_x}^\ast$, and take $\gamma_x = \beta^{-1}_x \alpha_x$.   

Since $\alpha \alpha^\ast$ is a natural transformation from $F$ to
itself, if we define $P(\alpha \alpha^\ast)_x = P(\alpha_x
{\alpha_x}^\ast)$ for any polynomial $P$, we have 
\[     P(\alpha_x {\alpha_x}^\ast) F(f) = F(f) P(\alpha_y {\alpha_y}^\ast)
\]
for any morphism $f \maps x \to y$.  By the finite-dimensional
spectral theorem, we can find a sequence of polynomials $P_i$ such
that $P_i(\alpha_x {\alpha_x}^\ast) \to \beta_x$ and $P_i(\alpha_y
{\alpha_y}^\ast) \to \beta_y$.  Thus
\[      \beta_x F(f) = F(f) \beta_y  , \]
so $\beta$ is a natural transformation from $F$ to itself.
It follows that $\gamma = \beta^{-1}\cdot \alpha$ is a natural transformation
from $F$ to $G$.   Clearly $\beta$ is self-adjoint and $\gamma$ is
unitary.  \qed

\subsection{Duality for morphisms}

Duals of morphisms in $2\Hilb$ are just adjoint functors.  
Normally one needs to distinguish between left and right adjoint
functors, but duality at the 2-morphism level allows us to turn left
adjoints into right adjoints, and vice versa:
\begin{prop}\et \label{2cat1} Suppose $F \maps H \to H'$, $G \maps H'
\to H$ are 
morphisms in $2\Hilb$.  Then $F$ is left adjoint to $G$ with unit
$\iota \maps 1_H \To FG$ and counit $\eps \maps GF \To 1_{H'}$  
if and only if $F$ is right adjoint to $G$ with unit $\eps^\ast \maps 1_{H'} 
\To GF$ and counit $\iota^\ast \maps FG \To 1_H$. \end{prop}
Proof - The triangle equations for $\iota$ and $\eps$:
\[     (\iota F) \cdot (F \eps) = 1_F, \qquad 
(G \iota) \cdot (\eps G) = 1_G ,\]
become equivalent to those for $\eps^\ast$ and $\iota^\ast$: 
\[  (\eps^\ast G) \cdot  (G \iota^\ast) = 1_G, \qquad 
(F \eps^\ast) \cdot (\iota^\ast F) = 1_F, \]
by taking duals.    \qed
\noindent As noted by Dolan \cite{Dolan}, it is probably quite generally
true in $n$-categories that duality for $j$-morphisms allows us to turn
`left duals' of $(j-1)$-morphisms into `right duals' and vice versa.  
This should give the theory of $n$-Hilbert spaces
quite a different flavor from general $n$-category theory.

Every morphism in $2\Hilb$ has an adjoint.   We prove this
using bases and the concept of a skeletal 2-Hilbert space.
\begin{defn}\et  A category is {\rm skeletal} if all isomorphic
objects are equal.  \end{defn}
\begin{defn}\et A {\rm unitary equivalence} between 2-Hilbert spaces
$H$ and $H'$ consists of morphisms $U \maps H \to H'$, $V \maps H' \to
H$ and unitary natural transformations $\iota \maps 1_H \To UV$, $\eps
\maps VU \To 1_{H'}$ forming an adjunction.  If there exists a unitary
equivalence between $H$ and $H'$, we say they are {\rm unitarily
equivalent}. \end{defn} 
\begin{prop}\et  \label{2cat2} Any 2-Hilbert space is 
unitarily equivalent to a skeletal 2-Hilbert space.  
\end{prop}
Proof - Let $\{e_\lambda\}$ be a basis for the 2-Hilbert space $H$.
For any nonnegative integers 
$\{n^\lambda\}$ with only finitely many nonzero, make a choice of direct sum
\[ \bigoplus_\lambda n^\lambda e_\lambda ,\]
where $n^\lambda e_\lambda$ denotes the direct sum of $n^\lambda$ copies of
$e_\lambda$.  (Recall that the direct sum is an object equipped with
particular morphisms; it is only unique up to isomorphism, but here we
fix a particular choice.)  Let $H_0$ denote the full subcategory of
$H$ with only these direct sums as objects.  Note that $H_0$ inherits
a 2-Hilbert space structure from $H$, and it is skeletal.   Let
$V \maps H_0 \to H$ denote the inclusion functor.  

For any $x \in H$ there is a unique object $U(x) \in H_0$ for which 
$V(U(x))$ is isomorphic to $x$.   
By Proposition \ref{H*4}, we may choose a unitary 
isomorphism 
\[             \iota_x \maps x \to V(U(x)). \] 
For $x = V(y)$ we have $U(x) = y$, so we choose $\iota_x$ to be
the identity in this case.
For each morphism $f\maps x
\to y$ define $U(f) \maps U(x) \to U(x')$ so that the following diagram
commutes: 
\[  \begin{diagram}[V(U(x))]
\node{x} \arrow{e,t}{f} \arrow{s,l}{\iota_x}  \node{y}
\arrow{s,r}{\iota_{y}} \\ 
\node{V(U(x))}  \arrow{e,b}{V(U(f))} \node{V(U(y))} 
\end{diagram} 
\]
It follows that $U \maps H \to H_0$ is a functor.  

One may check that $U$ and $V$ are actually morphisms of 2-Hilbert
spaces.  Moreover, one may check that there is a natural isomorphism
\[          \hom(Ux,y) \iso \hom(x,Vy)   \]
given by 
\[             f \mapsto \iota_x V(f). \]
It follows that $U$ is left adjoint to $V$. The unit of this
adjunction is $\iota$, while the counit is the identity.
These are both unitary natural transformations. \qed

Just as with Hilbert spaces, phrasing definitions and theorems about
2-Hilbert spaces in terms of a basis is usually a mistake, since
they should be manifestly invariant under unitary equivalence.  In
comparison, the use of bases to prove theorems is at worst a minor
lapse of taste, and sometimes convenient.  This is facilitated by the
use of skeletal 2-Hilbert spaces.

\begin{prop}\et  \label{2cat3} Let $F \maps H \to H'$ be a morphism in
$2\Hilb$.  Then there is a morphism $F^\ast \maps H' \to H$ that is left
and right adjoint to $F$.    \end{prop}
Proof - Here we opt for a lowbrow proof using bases, to illustrate the
analogy between an adjoint functor and the adjoint of a matrix.  
By Proposition \ref{2cat2} it suffices to
consider the case where $H$ and $H'$ are skeletal.  Let 
$\{e_\lambda\}$ be a basis for $H$ and $\{e'_\mu\}$ a basis for $H'$.  
Write 
\[             F(e_\lambda) = \bigoplus_\mu F_{\lambda \mu} e'_\mu
\]
where $F_{\lambda\mu}$ are nonnegative integers and $F_{\lambda\mu}
e'_\mu$ denotes the direct sum of $F_{\lambda\mu}$ copies of
$e'_\mu$.  Let 
\[             {F_{\mu\lambda}}^\ast = F_{\lambda\mu}. \]
Defining
\[           F^\ast(e'_\mu) = \bigoplus_{\mu} F^\ast_{\mu\lambda}
e_\lambda ,\]
one may check that $F^\ast$ extends uniquely to a morphism from $H'$
to $H$.   Note that both $\hom(F e_\lambda,
e'_\mu)$ and $\hom(e_\lambda, F^\ast e'_\mu)$ may be naturally
identified with a direct sum of $F_{\lambda\mu}$ copies of $\C$, which
sets up an isomorphism  $\hom(F e_\lambda, e'_\mu) \iso \hom(e_\lambda,
F^\ast e'_\mu)$.  One can check that this extends uniquely to a natural
isomorphism 
\[            \hom(Fx,y) \iso \hom(x,F^\ast y), \]
so $F^\ast$ is a right adjoint, and by Proposition \ref{2cat1} also a left
adjoint, of $F$.  \qed

A basic fact in Hilbert space theory is that two objects in $\Hilb$ are
isomorphic if and only if there is a unitary morphism between them.  The
same is true of objects in any other 2-Hilbert space, by Proposition
\ref{H*4}.  Similarly, two morphisms in $2\Hilb$ are isomorphic if and
only if there is a unitary natural transformation between them, by
Proposition \ref{unitary.nat}.  Below we show a similar result for
objects in $2\Hilb$.  In general, we expect a recursively defined
notion of `equivalence' of $j$-morphisms in an $n$-category: two
$n$-morphisms are equivalent if they are equal, while two
$(j-1)$-morphisms $x,y$ are equivalent if there exist $f \maps x \to
y$ and $g \maps y \to x$ with $gf$ and $fg$ equivalent to the identity
on $x$ and $y$, respectively.  In an $n$-Hilbert space we also expect
a similar notion of `unitary equivalence': two $(n-1)$-morphisms are
unitarily equivalent if they are equal, while two $(j-1)$-morphisms
$x,y$ are unitarily equivalent if there exists $u \maps x \to y$ with
$uu^\ast$ and $u^\ast u$ unitarily equivalent to $1_x$ and $1_y$,
respectively.  Our results so far lead us to suspect that, quite
generally, equivalent $j$-morphisms in an $n$-Hilbert space will be
unitarily equivalent.

\begin{defn}\label{equivalence}\et An {\rm equivalence} between
2-Hilbert spaces $H$ and $H'$ is an pair of morphisms $F \maps H \to
H'$, $G \maps H' \to H$ together with natural isomorphisms $\alpha \maps
1_H \To FG$, $\beta \maps GF \To 1_{H'}$.  If there is an equivalence
between $H$ and $H'$, we say they are {\rm equivalent}. \end{defn}

Note that a unitary equivalence is automatically an equivalence.
Conversely:

\begin{prop} \label{unitary.eq} \et Suppose $H$ and $H'$ are 2-Hilbert
spaces and the morphisms $F \maps H \to H'$, $G \maps H' \to H$ can be
extended to an equivalence between $H$ and $H'$.  Then $F$ and $G$ can
be extended to a unitary equivalence between $H$ and $H'$.
\end{prop}

Proof - Suppose $\alpha \maps 1_H \To FG$, $\beta \maps GF \To 1_{H'}$
are natural isomorphisms.  By Proposition \ref{unitary.nat} we can find
unitary natural transformations $\gamma \maps 1_H \To FG$, $\delta
\maps GF \To 1_{H'}$.   We may then obtain 
an adjunction by replacing $\gamma$ with the composite $\gamma'$ given by
\[  \begin{diagram} [FG = F1_{H'}G]
\node{1_H} \arrow{e,t}{\gamma} 
\node{FG = F1_{H'}G} \arrow{e,t}{F\delta^{-1}G} 
\node{FGFG} \arrow{e,t}{\gamma^{-1}FG} 
\node{FG} 
\end{diagram} \]
Checking that this is an adjunction is a lengthy but straightforward
calculation.  Noting that $\gamma'$ is unitary, we conclude that
$(F,G,\iota,\epsilon)$ is a unitary equivalence.  \qed

\noindent When we are being less pedantic, we call a 2-Hilbert
space morphism $F \maps H \to H'$ an {\it equivalence} if it can
be extended to an equivalence in the sense of Definition
\ref{equivalence}.  

Just as Hilbert spaces are classified by their dimension, we have:

\begin{cor} \label{2cat4} \et Two 2-Hilbert spaces are equivalent
if and only if they have the same dimension.
\end{cor}

Proof - Since an equivalence between $H$ and $H'$ carries a basis
of $H$ to a basis of $H'$, Proposition \ref{2hilb2} implies that
dimension is preserved by equivalence.  By Proposition
\ref{2cat2} it thus suffices to show two skeletal 2-Hilbert spaces are
equivalent if they have the same dimension.  Let
$\{e_\lambda\}$ be a basis of $H$ and $\{e'_\lambda\}$ a corresponding
basis of $H'$.  Then there is a unique 2-Hilbert space morphism with
$F(e_\lambda) = e'_\lambda$, and the adjunction constructed as in the
proof of Proposition \ref{2cat3} is a unitary equivalence. \qed

\subsection{Duality for objects} \label{dualityforobjects}

Finally, duals of objects in $2\Hilb$ are defined using an `internal
hom'.  Given 2-Hilbert spaces $H$ and $H'$, let $\hom(H,H')$ be the
category having 2-Hilbert space morphisms $F\maps H \to H'$ as objects
and 2-morphisms between these as morphisms.

\begin{prop} \label{2cat5}\et Suppose $H$ is a finite-dimensional
2-Hilbert space and $H'$ is a 2-Hilbert space.  Then the category
$\hom(H,H')$ becomes a
$\Hilb$-category if for any $F,G \in \hom(H,H')$ 
we make $\hom(F,G)$ into a Hilbert space with the
obvious linear structure and the inner product given by
\[           \langle \alpha, \beta \rangle = \sum_{\lambda} 
\langle \alpha_{e_\lambda} , \beta_{e_\lambda} \rangle \]
for any basis $\{e_\lambda\}$ of $H$. 
Moreover, $\hom(H,H')$ becomes a 2-Hilbert space if we define the dual
of $\alpha \maps F \To G$ by $(\alpha^\ast)_x = (\alpha_x)^\ast$. \end{prop}

Proof - Note first that $\hom(F,G)$ becomes a vector space if we define
\[        (\alpha + \beta)_x = \alpha_x + \beta_x, \qquad
            (c\alpha)_x = c(\alpha_x)  \]
for any $\alpha, \beta \maps F \To G$ and $c \in \C$.  Note also
that the inner product described above is nondegenerate, since if
$\alpha_{e_\lambda} = 0$ for all objects $e_\lambda$ in a basis,
then $\alpha = 0$.  Finally, note that the inner product is
independent of the choice of basis: if $\{e'_\lambda\}$ is
another basis we may assume after reordering that $e_\lambda \iso
e'_\lambda$, and by Proposition \ref{H*4} we may choose unitary
isomorphisms $u_\lambda \maps e_\lambda \to e'_\lambda$, so that
\[ \alpha_{e'_\lambda} = F(u_\lambda)^\ast \alpha_{e_\lambda}
G(u_\lambda) \]
and similarly for $\beta$.  It follows that
\ban  \langle \alpha_{e'_\lambda} , \beta_{e'_\lambda} \rangle &=&
\langle F(u_\lambda)^\ast \alpha_{e_\lambda} G(u_\lambda) , 
F(u_\lambda)^\ast \beta_{e_\lambda} G(u_\lambda) \rangle \\
&=& \langle \alpha_{e_\lambda} , \beta_{e_\lambda} \rangle .\ean
Since composition of morphisms in $\hom(H,H')$ is bilinear, it becomes a
$\Hilb$-category.  

It is easy to check that defining $(\alpha^\ast)_x =
(\alpha_x)^\ast$ makes $\hom(H,H')$ into a $\ast$-category, and
using Proposition \ref{H*1} one can also check it is an H*-category. 
To check that it is a 2-Hilbert space it suffices by Proposition
\ref{2hilb1} to check that it has a zero object, direct sums and
kernels.  Any functor $0 \maps H \to H'$ mapping all objects in
$H$ to zero objects in $H'$ is initial in  $\hom(H,H')$.  Given
$F,F' \in \hom(H,H')$, we may take as the direct sum $F \oplus
F'$ any functor with $(F \oplus F')(x) = F(x) \oplus F(x')$ for
any object $x \in H$ and $(F \oplus F')(f) = F(f) \oplus F(f')$
for any morphism $f$.  Similarly, given $\alpha \maps F \to F'$,
we may construct $\ker \alpha \in \hom(H,H')$ by letting $(\ker
\alpha)(x) = \ker \alpha_x$ for any object $x$ and defining
$(\ker \alpha)(f)$ for any morphism using the universal property
of the kernel.  \hbox{\hskip 30em} \qed

\begin{defn} \et Given a finite-dimensional 2-Hilbert space $H$, the
{\rm dual} $H^\ast$ is the 2-Hilbert space $\hom(H,\Hilb)$.  \end{defn}

The following is an analog of the Riesz representation theorem for
finite-dimensional 2-Hilbert spaces.  In its finite-dimensional form,
the Riesz representation theorem says if $x$ is a Hilbert space, any
morphism $f \maps x \to \C$ is equal to one of the form 
\[                 \langle v, \cdot  \rangle \]
for some $v \in H$.  This determines an isomorphism $\overline x \iso
x^\ast$.  Similarly, given a 2-Hilbert space $H$, we say a
morphism $F \maps H \to \Hilb$ is {\it representable} if it is naturally
isomorphic to one of the form 
\[                 \hom(x, \cdot)  \]
for some $x \in H$.   The essence of the Riesz
representation theorem for 2-Hilbert spaces is that every morphism
$F \maps H \to \Hilb$ is representable.  This yields an 
equivalence between $H^\op$ and $H^\ast$.  

\begin{prop} \et  For any finite-dimensional 2-Hilbert space $H$, the
morphism \break $U \maps H^\op \to H^\ast$ given by 
\[        U(x) = \hom(x,\cdot) , \qquad  U(f) = \hom(f,\cdot)  \]
is an equivalence between $H^\op$ and $H^\ast$.
\end{prop} 

Proof - It suffices to show that $U$ is fully faithful and
essentially surjective.  We can check both of these using a
basis $\{e_\lambda\}$ of $H$.   We leave the full faithfulness to the
reader.  Checking that $U$ is essentially surjective amounts to checking
that any $F \in H^\ast$ is representable.  Note there is a `dual
basis' of 2-Hilbert space morphisms $f^\lambda \in \hom(H,\Hilb)$ with
\[ f^\mu(e_\lambda) \iso \cases{\C & $\lambda = \mu$ \cr
                     0  & $\lambda \ne \mu$ } \]
Since any morphism $F \maps H \to \Hilb$ is determined up
to natural isomorphism by its value on the basis $\{e_\lambda\}$, any $F \in
H^\ast$ is isomorphic to a direct sum of the $\{f^\lambda\}$.  But 
$f^\lambda$ is isomorphic to $U(e_\lambda)$, so $U$ is essentially
surjective.  \qed

\subsection{The tensor product}\label{tensorproduct}

Next we develop the tensor product of 2-Hilbert spaces.  For this we need
the analog of a bilinear map:
\begin{defn} \et Given 2-Hilbert spaces $H,H',K$, a functor
$F \maps H \times H' \to K$ is a {\rm bimorphism} of 2-Hilbert spaces
if for any objects $x \in H$, $x' \in H'$ the functors $F(x \tensor
\,\cdot\,) \maps H' \to K$ and $F(\,\cdot \,\tensor x') \maps H \to K$
are 2-Hilbert space morphisms.   We write $\bihom(H\times H', K)$ for
the category having bimorphisms $F \maps H \times H' \to K$ as objects
and natural transformations between these as morphisms. 
\end{defn}  

\begin{prop} \et Suppose $H$ and $H'$ are finite-dimensional
2-Hilbert spaces and $K$ is a 2-Hilbert space.  Then $\bihom(H\times H',K)$ 
becomes a $\Hilb$-category if for any $F,G \in \bihom(H \times H',K)$ 
we make $\hom(F,G)$ into a Hilbert space with the
obvious linear structure and the inner product given by
\[           \langle \alpha, \beta \rangle = \sum_{\lambda,\mu} 
\langle \alpha_{(e_\lambda,f_\mu)} , \beta_{(e_\lambda,f_\mu)} \rangle \]
for any bases $\{e_\lambda\}$ of $H$ and $\{f_\mu\}$ of $H'$. 
Moreover, $\bihom(H \times H',K)$ becomes a 2-Hilbert space if we define
the dual of $\alpha \maps F \To G$ by $(\alpha^\ast)_x =
(\alpha_x)^\ast$. \end{prop} 

Proof - The proof is analogous to that of Proposition \ref{2cat5}.  \qed

Given 2-Hilbert spaces $H,H'$ and $L$, note that a bimorphism
$T \maps H \times H' \to L$ induces a morphism 
\[            T^\ast \maps \hom(L,K) \to \bihom(H \times H',K). \]
\begin{defn} \et  Given 2-Hilbert space $H,H'$, a {\rm tensor product}
of $H$ and $H'$ is a bimorphism $T \maps H \times H' \to L$ together
with a choice for each 2-Hilbert space $K$
of an equivalence of 2-Hilbert spaces extending $T^\ast \maps
\hom(L,K) \to \bihom(H \times H',K)$.   \end{defn}

\noindent In the above situation, by abuse of language we
may say simply that $T \maps H \times H'
\to L$ is a tensor product of $H$ and $H'$.

\begin{prop} \label{tensor} \et  Given finite-dimensional
2-Hilbert spaces $H$ and $H'$, there
exists a tensor product $T \maps H \times H' \to L$.  Given another
tensor product $T' \maps H \times H' \to L'$, there is an equivalence
$F \maps L \to L'$ for which the following diagram commutes up to
a specified natural isomorphism:
\[
\begin{diagram}[H \times H']
\node[2]{H \times H'} \arrow{sw,t}{T} \arrow{se,t}{T'} \\
\node{L} \arrow[2]{e,b}{F} \node[2]{L'}  
\end{diagram}
\]
\end{prop}
  
Proof - Let $\{e_\lambda\}$ be a basis for $H$, and
$\{f_\mu\}$ a basis for $H'$.  Let $L$ be the skeletal
2-Hilbert space with a basis of objects denoted by $\{e_\lambda \tensor
f_\mu\}$, and with
\[       \hom(e_\lambda \tensor f_\mu,e_\lambda \tensor f_\mu) = 
\hom(e_\lambda,e_\lambda) \tensor \hom(f_\mu,f_\mu)  \]
as H*-algebras (using the obvious tensor product of H*-algebras).
There is a unique bimorphism $T \maps H \times H' \to L$ with
$T(e_\lambda,f_\mu) = e_\lambda \tensor f_\mu$.  Given a 2-Hilbert
space $K$ one may check that $T^\ast \maps \hom(L,K) \to \bihom(H
\times H',K)$ extends to an equivalence.  Choosing such an equivalence
for every $K$ we obtain a tensor product of $H$ and $H'$.

Given two tensor products as in the statement of the proposition, let
$F\maps L \to L'$ be the image of $T'$ under the chosen equivalence
$\bihom(H \times H', L') \simeq \hom(L,L')$.  One can check that $L$
is an equivalence and that the above diagram commutes up to a
specified natural isomorphism, much as in the usual proof that the tensor
product of vector spaces is unique up to a specified isomorphism.
\qed

Given a tensor product of the 2-Hilbert spaces $H$ and $H'$, we often
write its underlying 2-Hilbert space as $H \tensor H'$.  This notation
may tempt one to speak of `the' tensor product of $H$ and $H'$, which is 
is legitimate if one uses the generalized 
`the' as advocated by Dolan \cite{Dolan}.  In a set, when we
speak of `the' element with a given property, we implicitly mean that
this element is unique.  In a category, when we speak of
`the' object with a given property, we merely mean that this object is
unique up to isomorphism --- typically a specified isomorphism.
Similarly, in a 2-category, when we speak of `the' object with a given
property, we mean that this object is unique up to equivalence ---
typically an equivalence that is specified up to a specified
isomorphism.  This is the sense in which we may refer to `the' tensor
product of $H$ and $H'$.  The generalized `the' may be extended in an
obvious recursive fashion to $n$-categories.

Suppose that $H$ and $H'$ are finite-dimensional 2-Hilbert spaces.
Then for any pair of objects $x \in H$, $x' \in H'$, 
we can use the bimorphism $T \maps H \times H' \to H \tensor H'$ 
to define an object $x \tensor x' = T(x,x')$ in $H \tensor H'$.
Similarly, given a morphism $f \maps x \to y$ in
$H$ and a morphism $f' \maps x' \to y'$, we obtain a morphism
\[   f \tensor f'\, \maps x \tensor x' \to y \tensor y'\]
in $H \tensor H'$.  We usually write
\[  f \tensor x' \, \maps x \tensor x' \to y \tensor x' \]
for the morphism $f \tensor 1_{x'}$, and 
\[  x \tensor f' \, \maps x \tensor x' \to x \tensor y' \]
for the morphism $1_x \tensor f'$.  

We expect that $2\Hilb$ has the structure of a monoidal 2-category with
the above tensor product as part of the monoidal structure.  Kapranov
and Voevodsky \cite{KV} have defined the notion of a weak monoidal
structure on a strict 2-category, which should be sufficient for the
purpose at hand.  On the other hand the work of Gordon, Power and Street
\cite{GPS} gives a fully general notion of weak monoidal 2-category,
namely a 1-object tricategory.  This should also be suitable for
studying the tensor product on $2\Hilb$, though it might be considered
overkill.  Both these sorts of monoidal 2-category involve various
extra structures besides the tensor product of objects in $2\Hilb$.
Most of these should arise from the universal property of the tensor
product.  

For example, suppose we are given a morphism $F \maps H \to H'$ and an
object $K$ in $2\Hilb$.    Thus we have bimorphisms $T \maps H \times
K \to H \tensor K$ and $T' \maps H \times K' \to H \tensor K'$, and 
$T^\ast$ has some morphism
\[         S \maps \bihom(H \times K, H' \tensor K) \to \hom(H \tensor
K,H' \tensor K) \]
as inverse up to natural isomorphism.  
Applying $S$ to the bimorphism given by the composite 
\[
\begin{diagram} [H \tensor K]
\node{H \times K} \arrow{e,t}{F \times 1_K} \node{H' \times K} \arrow{e,t}{T'}
\node{H' \tensor K} 
\end{diagram}
\]
we obtain a morphism we denote by
\[    F \tensor K \maps H \tensor K \to H' \tensor K .\]
Similarly, given an object $H \in 2\Hilb$ and a morphism $G \maps K \to
K'$, we obtain a morphism 
\[     H \tensor G \maps H \tensor K \to H \tensor K'. \]
Moreover, we have:

\begin{prop} \label{tensorator} \et Let $F \maps H \to H'$ and $G \maps K
\to K'$ be morphisms in $2\Hilb$.  Then the following diagram 
\[
\begin{diagram}
\node{H \tensor K} 
\arrow{e,t}{F \tensor K} \arrow{s,l}{H \tensor G}
\node{H' \tensor K} 
\arrow{s,r}{H' \tensor G} \\
\node{H \tensor K'} \arrow{e,b}{F \tensor K'} \node{H' \tensor K'} 
\end{diagram}
\]
commutes up to a specified natural isomorphism
\[  {\bigotimes}_{F,G} \maps (F \tensor K)(H' \tensor G) \To
(H \tensor G)(F \tensor K') . \]
\end{prop}

Proof - Here we have fixed tensor products of all the 2-Hilbert spaces
involved, so we have bimorphisms
\[        T_{H,K} \maps H \times K \to H \tensor K \]
and so on.   Applying the equivalence
\[          \bihom(H \times K, H' \tensor K') \simeq \hom(H \tensor
K,H' \tensor K') \]
coming from the definition of tensor product to the bimorphism given by
the composite  
\be
\begin{diagram} [H \tensor K]
\node{H \times K} \arrow{e,t}{F \times G} \node{H' \times K'}
\arrow{e,t}{T_{H',K'}} \node{H' \tensor K'} 
\end{diagram} \label{bigotimes1}
\ee
we obtain a morphism we denote by
\[         F \tensor G \maps H \tensor K \to H' \tensor K'.  \]
We shall construct a natural isomorphism from $(F \tensor K)(H' \tensor G)$ 
to $F \tensor G$.  Composing this with an analogous natural isomorphism
from $F \tensor G$ to $(H \tensor G)(F \tensor K')$ one obtains
$\bigotimes_{F,G}$.  

If we precompose $F \tensor G$ with $T_{H,K}$ we obtain a bimorphism
naturally isomorphic to (\ref{bigotimes1}).  
If we precompose $(F \tensor K')(H \tensor G)$ with $T_{H,K}$, we obtain
a bimorphism naturally isomorphic to 
\be
\begin{diagram} [H' \tensor K'.]
\node{H \times K} \arrow{e,t}{F \times K} \node{H' \times K}
\arrow{e,t}{T_{H',K}} \node{H' \tensor K}  
\arrow{e,t}{H' \tensor G} \node{H' \tensor K'}
\end{diagram} \label{bigotimes2}
\ee
Note also that in both cases, a {\it specified}
natural isomorphism is given by the definition of tensor product.
Since precomposition with $T_{H,K}$ is an equivalence between 
$\bihom(H \times K,H' \tensor K')$ and $\hom(H \tensor K,H' \tensor K'),$
it thus suffices to exhibit a natural isomorphism between
(\ref{bigotimes1}) and (\ref{bigotimes2}).

Factoring these by $F \times K$, it suffices to exhibit a natural
isomorphism between 
\[
\begin{diagram} [H' \tensor K]
\node{H' \times K} \arrow{e,t}{H' \times G} 
\node{H' \tensor K'} \arrow{e,t}{T_{H',K'}}
\node{H' \tensor K'} 
\end{diagram} 
\]
and
\[
\begin{diagram} [H' \tensor K'.]
\node{H' \times K}
\arrow{e,t}{T_{H',K}} \node{H' \tensor K}  
\arrow{e,t}{H' \tensor G} \node{H' \tensor K'}
\end{diagram} 
\]
This arises from the definition of $H' \tensor G$.   \qed

\noindent The 2-morphism $\bigotimes_{F,G}$ is part of the
structure one  expects in a monoidal 2-category, and the fact
that the diagram in Proposition \ref{tensorator} does not commute
`on the nose' is one of the key ways in which monoidal
2-categories differ from monoidal categories.

We expect a 2-categorical version of $\hom$-tensor adjointness to hold
for the tensor product defined in this section and the $\hom$ defined
in section \ref{dualityforobjects}.  In other words, given
finite-dimensional 2-Hilbert space $H, H',$ and $K$, the obvious
functor from $\hom(H, \hom(H',K))$ to $\hom(H \tensor H', K)$ should
be an equivalence.  However, we shall not prove this here.  

\subsection{The braiding}

The symmetry in $\Cat$ gives braiding morphisms in $2\Hilb$ as follows. 
Let $H$ and $H'$ be 2-Hilbert spaces.  We may take their tensor product
in either order, obtaining tensor products $T \maps H \times H' \to H
\tensor H'$ and $T' \maps H' \times H \to H' \tensor H$.  
By the universal property of the tensor product, the bimorphism given by
the composite
\[
\begin{diagram} [H \tensor K]
\node{H \times H'} \arrow{e,t}{S_{H,H'}} \node{H' \times H} \arrow{e,t}{T'}
\node{H' \tensor H} 
\end{diagram}
\]
defines a morphism, the {\it braiding} 
\[   R_{H,H'} \maps H \tensor H' \to H' \tensor H .\]
One can check that $R_{H,H'}$ is an equivalence.  

We expect that $2\Hilb$ has the structure of a braided monoidal
2-category with the above braiding morphisms. However, the existing
notion of semistrict braided monoidal 2-category introduced by Kapranov
and Voevodsky \cite{KV} and subsequently refined in HDA1 is
insufficiently general to cover this example, since $2\Hilb$ is not a
semistrict monoidal 2-category.  One should however be able to
strictify $2\Hilb$, obtaining a semistrict braided monoidal
2-category.  Alternatively, the work of Trimble \cite{Trimble} should
give a fully general notion of weak braided monoidal 2-category, namely
a tetracategory with one object and one morphism.  This should apply to
$2\Hilb$ without strictification.

In any event, both semistrict and weak braided monoidal 2-categories
involve various structures in addition to the braiding morphisms.  Most
of these should arise from the universal property of the tensor product
together with the properties of the symmetry in $\Cat$.  For example, we
have:

\begin{prop} \label{braiding.naturalizer} \et Let $F \maps H \to H'$ be a
morphism and let $K$ be an object in $2\Hilb$.  Then the following
diagram 
\[
\begin{diagram}
\node{H \tensor K} 
\arrow{e,t}{F \tensor K} \arrow{s,l}{R_{H,K}}
\node{H' \tensor K} 
\arrow{s,r}{R_{H',K}} \\
\node{K \tensor H} \arrow{e,b}{K \tensor F} \node{K \tensor H'} 
\end{diagram}
\]
commutes up to a specified natural isomorphism 
\[      R_{F,K} \maps (F \tensor K) R_{H',K} \To R_{H,K} (K \tensor F) .\]
Similarly, given an object $H$ and a morphism $G \maps K
\to K'$ in $2\Hilb$, the following diagram 
\[
\begin{diagram}
\node{H \tensor K} 
\arrow{e,t}{H \tensor G} \arrow{s,l}{R_{H,K}}
\node{H \tensor K'} 
\arrow{s,r}{R_{H,K'}} \\
\node{K \tensor H} \arrow{e,b}{G \tensor H} \node{K' \tensor H} 
\end{diagram}
\]
commutes up to a specified natural isomorphism
\[      R_{H,G} \maps (H \tensor G)R_{H,K'} \To R_{H,K}(G \tensor H)
.\]
\end{prop}

Proof - We only treat the first case as the second is analogous. 
Applying the equivalence
\[         \bihom(H \times K, K' \tensor H') \simeq \hom(H \tensor
K,K' \tensor H') \]
coming from the definition of tensor product to the bimorphism given by
the composite  
\be
\begin{diagram} [H \tensor K]
\node{H \times K} \arrow{e,t}{F \times K} 
\node{H' \times K} \arrow{e,t}{S_{H',K}}
\node{K \times H'} \arrow{e,t}{T_{K,H'}} 
\node{K \tensor H'}   
\end{diagram} \label{braiding.nat.1}
\ee
we obtain a morphism we denote by $A \maps H \tensor K \to K' \tensor H$.
We shall construct a natural isomorphism from $(F \tensor K)R_{H',K}$
to $A$. 
Using the fact that (\ref{braiding.nat.1}) equals
\[
\begin{diagram} [H \tensor K]
\node{H \times K} \arrow{e,t}{S_{H,K}} 
\node{K \times H} \arrow{e,t}{K \times F}
\node{K \times H'} \arrow{e,t}{T_{K,H'}} 
\node{K \tensor H'}   
\end{diagram} 
\]
one can similarly obtain a natural isomorphism from $A$ to
$R_{H,K}(K \tensor F)$.  The composite of these is $R_{F,K}$.  

If we precompose $A$ with $T_{H,K}$ we obtain a bimorphism
naturally isomorphic to (\ref{braiding.nat.1}).  
If we precompose $(F \tensor K)R_{H',K}$ with $T_{H,K}$, we obtain
a bimorphism naturally isomorphic to 
\be
\begin{diagram} [H' \tensor K.]
\node{H \times K} \arrow{e,t}{F \times K} \node{H' \times K}
\arrow{e,t}{T_{H',K}} \node{H' \tensor K}  
\arrow{e,t}{R_{H',K}} \node{K' \tensor H}
\end{diagram} \label{braiding.nat.2}
\ee
In both cases, a
natural isomorphism is given by the definition of tensor product.
It thus suffices to exhibit a natural isomorphism between
(\ref{braiding.nat.1}) and (\ref{braiding.nat.2}).  This may be
constructed as in the proof of Proposition \ref{tensorator}. \qed

\subsection{The involutor}

As indicated in Figure 1, for $2\Hilb$ to be a stable 2-category it
should possess an extra layer of structure after the tensor product
and the braiding, namely the `involutor'.  Also, this structure should
have an extra property making $2\Hilb$ `strongly involutory'.  The
involutor is a weakened form of the equation appearing in the
definition of a symmetric monoidal category.  Namely, while the
braiding need not satisfy
\[             R_{H',H} R_{H,H'} = 1_{H \tensor H'} \]
for all objects $H,H' \in 2\Hilb$, there should be a 2-isomorphism
\[          I_{H,H'} \maps R_{H,H'} R_{H',H} \To 1_{H \tensor H'}, \]
the involutor.   

We construct the involutor as follows.  Choose tensor
products $T \maps H \times H' \to H \tensor H'$ and $T' \maps H' \times
H \to H' \tensor H$. 
Then by the universality of the tensor product, the
commutativity of 
\[
\begin{diagram} [H \times H']
\node[2]{H' \times H} \arrow{se,t}{S_{H',H}}  \\
\node{H \times H'} \arrow[2]{e,b}{1_{H \times H'}}
\arrow{ne,t}{S_{H,H'}} \node[2]{H \times H'}   \\
\end{diagram}
\]
implies that 
\[
\begin{diagram} [H \times H']
\node[2]{H' \otimes H} \arrow{se,t}{R_{H',H}}  \\
\node{H \otimes H'} \arrow[2]{e,b}{1_{H \otimes H'}}
\arrow{ne,t}{R_{H,H'}} \node[2]{H \otimes H'}   \\
\end{diagram}
\]
commutes up to a specified natural transformation.  This is the involutor
\[          I_{H,H'} \maps R_{H,H'} R_{H',H} \To 1_{H \tensor H'}. \]

In addition, for $2\Hilb$ to be stable, or `strongly involutory', the
involutor should satisfy a special coherence law of its own, in
analogy to how the braiding satisfies a special equation in a
symmetric monoidal category.  In HDA0 this equation was described in
terms of $R_{H,H'}$ and a weak inverse thereof, but it turns out
to be easier to give the equation by stating that
the following horizontal composites agree:
\[ I_{H,H'} \circ 1_{R_{H,H'}} \maps R_{H,H'} R_{H',H} R_{H,H'}
\To R_{H,H'} \]
and 
\[  1_{R_{H,H'}} \circ I_{H,H'} \maps R_{H,H'} R_{H',H} R_{H,H'}
\To R_{H,H'} \]
This is indeed the case, as one can show using the properties of
the tensor product.

\section{2-H*-algebras}\label{2-H*-algebras}

Now we consider 2-Hilbert spaces with extra structure and properties,
as listed in the second column of Figure 2.  

\begin{defn} \et  A {\rm 2-H*-algebra} $H$ is a 2-Hilbert space
equipped with a {\rm product} bimorphism $\tensor \maps H \times
H \to H$, a {\rm unit} object $1 \in H$, a unitary natural
transformation $a_{x,y,z} \maps (x \tensor y) \tensor z \to x
\tensor (y \tensor z)$ called the {\rm associator}, and unitary
natural transformations $\ell_x \maps 1 \tensor x \to x$, $r_x
\maps x \tensor 1 \to x$ called the {\rm left} and {\rm right
unit laws}, making $H$ into a monoidal category.  We require also
that every object $x \in H$ has a left dual.  \end{defn}

\noindent Recall that for $H$ to be a monoidal category, one demands that
the following pentagon commute:
\[    \begin{diagram}[((x \tensor y)\tensor z)\tensor w]
\node{((x \tensor y) \tensor z)\tensor w}
\arrow{e,t}{a_{x\tensor y,z,w}}
\arrow{s,l}{a_{x,y,z}\tensor w}
\node{(x\tensor y)\tensor(z\tensor w)}
\arrow{e,t}{a_{x,y,z\tensor w}}
\node{x \tensor (y\tensor(z \tensor w))}  \\
\node{(x \tensor (y \tensor z)) \tensor w}
\arrow[2]{e,t}{a_{x,y\tensor z,w}}
\node[2]{x \tensor ((y \tensor z)\tensor w)}
\arrow{n,r}{x \tensor a_{y,z,w}}
\end{diagram} \]
as well as the following diagram involving the unit laws:
\[      \begin{diagram}[(1 \tensor x) \tensor 1]
\node{(1 \tensor x) \tensor 1} \arrow{s,l}{\ell_x \tensor 1} 
\arrow[2]{e,t}{a_{1,x,1}} \node[2]{1 \tensor (x\tensor 1)} \arrow{s,r}{1
\tensor r_x} \\
\node{x \tensor 1} \arrow{e,t}{r_x} \node{x} \node{1 \tensor x}
\arrow{w,t}{\ell_x} 
\end{diagram} \]

Mac Lane's coherence theorem \cite{MacLane2} says that every monoidal
category is equivalent, as a monoidal category, to a {\it strict}
monoidal category, that is, one for which the associators and unit
laws are all identity morphisms.  Sometimes we will use this to
streamline formulas by not parenthesizing tensor products and not
writing the associators and unit laws.  Such formulas apply literally
only to the strict case, but one can always use Mac Lane's theorem to
apply them to general monoidal categories.  In practice, this amounts
to parenthesizing tensor products however one likes, and inserting 
associators and unit laws when needed to make the formulas make sense.  

A {\it left dual} of an object $x$ in a monoidal category
is an object $y$ together with morphisms 
\[          e \maps y \tensor x \to 1 \]
and 
\[           i \maps 1 \to x \tensor y ,\]
called the {\it unit} and {\it counit}, such that the following 
diagrams commute:
\[
\begin{diagram} [x \tensor y \tensor x]
\node{x} \arrow[2]{e,t}{1_x} \arrow{se,b}{i \tensor x}
\node[2]{x}   \\
\node[2]{x \tensor y \tensor x} \arrow{ne,b}{x \tensor e}
\end{diagram}
\]
\[
\begin{diagram} [x \tensor y \tensor x]
\node{y} \arrow[2]{e,t}{1_y} \arrow{se,b}{y \tensor i}
\node[2]{y}   \\
\node[2]{y \tensor x \tensor y} \arrow{ne,b}{e \tensor y}
\end{diagram}
\]
(These diagrams apply literally only when the monoidal category
is strict.)  In this situation we also say that $x$ is a {\it right dual}
of $y$, and that $(x,y,i,e)$ is an {\it adjunction}.  All adjunctions
having $x$ as right dual are uniquely isomorphic in the following sense:

\begin{prop} \label{2H*1} \et  Given an adjunction $(x,y,i,e)$
in a monoidal category and an isomorphism $f\maps y \to y'$, there is
an adjunction $(x,y',i',e')$ given by:
\[       i' = i(x \tensor f), \qquad e' = (f^{-1} \tensor x)e.\]
Conversely, given two adjunctions $(x,y,i,e)$ and $(x,y',i',e')$,
there is a unique isomorphism $f \maps y \to y'$
for which $i' = i(x \tensor f)$ and $e' = (f^{-1} \tensor x)e$.
This is given in the strict case by the composite
\[  \begin{diagram}[y \tensor (x \tensor y')]
\node{y = y \tensor 1} \arrow{e,t}{y \tensor i'}
\node{y \tensor x \tensor y'} \arrow{e,t}{e \tensor y'}
\node{1 \tensor y' = y'}
\end{diagram} \]
\end{prop}

Proof - This result is well-known and the proof is a simple calculation.
\qed 

\noindent  Similarly, any two adjunctions having a given object as
right dual are canonically isomorphic.  We may thus
speak of `the' left or right dual of a given object, using the 
generalized `the', as described in Section \ref{tensorproduct}.
Note that duality at the morphism level of a 2-H*-algebra allows us to turn 
left duals into right duals, and vice versa, at the object level: 

\begin{prop}\label{2H*2}\et Suppose that $H$ is a 2-H*-algebra.  Then 
$(x,x^\ast,i,e)$ is an adjunction if and only if $(x^\ast,x,e^\ast,i^\ast)$
is an adjunction.  \end{prop}

Proof - The proof is analogous to that of Proposition \ref{2cat1}.  \qed

Next we turn to braided and symmetric 2-H*-algebras.
A good example of a braided 2-H*-algebra is the category
of tilting modules of a quantum group when the parameter $q$ is a
suitable root of unity \cite{CP}.  Categories very similar to our
braided 2-H*-algebras have been studied by Fr\"ohlich and Kerler
\cite{FK} under the name `C*-quantum categories'; our definitions differ
only in some fine points.  A good example of a symmetric 2-H*-algebra
is the category of finite-dimensional continuous unitary
representations of a compact topological group.  Doplicher and Roberts
\cite{DR} have studied categories very similar to our symmetric
2-H*-algebras.

\begin{defn} \et A {\rm braided 2-H*-algebra} is a 2-H*-algebra $H$
equipped with a unitary natural isomorphism $B_{x,y} \maps x \tensor y
\to y \tensor x$ making $H$ into a braided monoidal category.
\end{defn}

\begin{defn} \et A {\rm symmetric 2-H*-algebra} is a 2-H*-algebra for
which the braiding is a symmetry.  \end{defn}

Recall that for $H$ to be a braided monoidal category, the 
following two hexagons must commute:
\[ \begin{diagram}[(x \tensor y) \tensor z)]
\node{x \tensor (y \tensor z)} \arrow{s,l}{B_{x,y \tensor z}}
\arrow{e,t}{a^{-1}_{x,y,z}}
\node{(x \tensor y) \tensor z} \arrow{e,t}{B_{x,y} \tensor z}
\node{(y \tensor x) \tensor z} \arrow{s,r}{a_{y,x,z}} \\
\node{(y \tensor z) \tensor x} \arrow{e,b}{a_{y,z,x}}
\node{y \tensor (z \tensor x)} \arrow{e,t}{y \tensor B_{x,z}} 
\node{y \tensor (x \tensor z)} 
\end{diagram} \]

\[ \begin{diagram}[(x \tensor y) \tensor z)]
\node{(x \tensor y) \tensor z} \arrow{s,l}{B_{x \tensor y, z}}
\arrow{r,t}{a_{x,y,z}} 
\node{x \tensor (y \tensor z)} \arrow{e,t}{x \tensor B_{y,z}}
\node{x \tensor (z \tensor y)} \arrow{s,r}{a_{x,z,y}^{-1}} \\
\node{z \tensor (x \tensor y)} \arrow{e,b}{a^{-1}_{z,x,y}}
\node{(z \tensor x) \tensor y}  \arrow{e,t}{B_{x,z} \tensor y} 
\node{(x \tensor z) \tensor y}
\end{diagram} \]
as well as the following diagrams:
\[ \begin{diagram}[1 \tensor x]
\node{1 \tensor x} \arrow{se,b}{\ell_x} \arrow[2]{e,t}{B_{1,x}} 
\node[2]{x \tensor 1} \arrow{sw,b}{r_x} \\
\node[2]{x} 
\end{diagram} \]

\[ \begin{diagram}[1 \tensor x]
\node{x \tensor 1} \arrow{se,b}{r_x} \arrow[2]{e,t}{B_{x,1}} 
\node[2]{1 \tensor x} \arrow{sw,b}{\ell_x} \\
\node[2]{x} 
\end{diagram} \]
The braiding is a symmetry if $B_{x,y} =
B_{y,x}^{-1}$ for all objects $x$ and $y$.

\subsection{The balancing} \label{balancing}

In the study of braided monoidal categories where objects have duals,
it is common to introduce something called the `balancing'.  The
balancing can treated in various ways \cite{FK,JS2,RT}.  For example, one
may think of it as a choice of automorphism $b_x \maps x \to x$ for
each object $x$, which is required to satisfy certain laws.
While very important in topology, this extra structure seems somewhat
ad hoc and mysterious from the algebraic point of view.  We now
show that braided 2-H*-algebras are automatically equipped with
a balancing.  The reason is that not only the objects, but also
the morphisms, have duals.  In fact, some of what follows would
apply to any braided monoidal category in which both objects and
morphisms have duals.

In any 2-H*-algebra, Proposition \ref{2H*2} gives a way to make
any object $x$ into the left dual of its left dual $x^\ast$.  In a braided
2-H*-algebra, $x$ also becomes the left dual of $x^\ast$ in another
way:

\begin{prop} \label{2H*2.5} \et 
Let $H$ be a braided 2-H*-algebra.  Then $(x,x^\ast,i,e)$
is an adjunction if and only if
$(x^\ast,x,iB_{x,x^\ast},B_{x,x^\ast}e)$ is an adjunction.  \end{prop}

Proof - The proof is a simple computation.  \qed

It follows from Proposition \ref{2H*1} that these two ways
to make $x$ into the left dual of $x^\ast$ determine an automorphism
of $x$.  Simplifying the formula for this automorphism
somewhat, we make the following definition:

\begin{defn}\et If $H$ is a braided 2-H*-algebra and
$(x,x^\ast,i,e)$ is an adjunction in $H$, the {\rm balancing} of the
adjunction is the morphism $b \maps x \to x$ given in the
strict case by the composite:
\[    \begin{diagram} [(x \tensor x^\ast) \tensor x]
\node{x} \arrow{e,t}{e^\ast \tensor x}
\node{x^\ast \tensor x \tensor x}
\arrow{e,t}{x^\ast \tensor B_{x,x}}
\node{x^\ast \tensor x \tensor x}
\arrow{e,t}{e \tensor x} \node{x} 
\end{diagram} \]
\end{defn}

It is perhaps easiest to understand the significance of the balancing
in terms of its relation to topology.  We shall be quite
sketchy about describing this, but the reader can fill in the details
using the ideas described in HDA0 and the many references therein.
Especially relevant is the work of Freyd and Yetter \cite{FY}, Joyal
and Street \cite{JS}, and Reshetikhin and Turaev \cite{RT,T}.  We
discuss this relationship more carefully in the Conclusions.

\bigskip
\centerline{\epsfysize=1.5in\epsfbox{2dtangle.eps}}
\medskip
\centerline{3.  Typical tangle in 2 dimensions}
\medskip

The basic idea is to use tangles to represent certain morphisms in
2-H*-algebras.  A typical oriented tangle in 2 dimensions is shown in
Figure 3.  If we fix an adjunction $(x,x^\ast,i,e)$ in a strict
2-H*-algebra $H$, any such tangle corresponds uniquely to a morphism
in $H$ as follows.  As shown in Figure 4, vertical juxtaposition of
tangles corresponds to the composition of morphisms, while horizontal
juxtaposition corresponds to the tensor product of morphisms.

\vbox{
\bigskip
\centerline{\epsfysize=1.5in\epsfbox{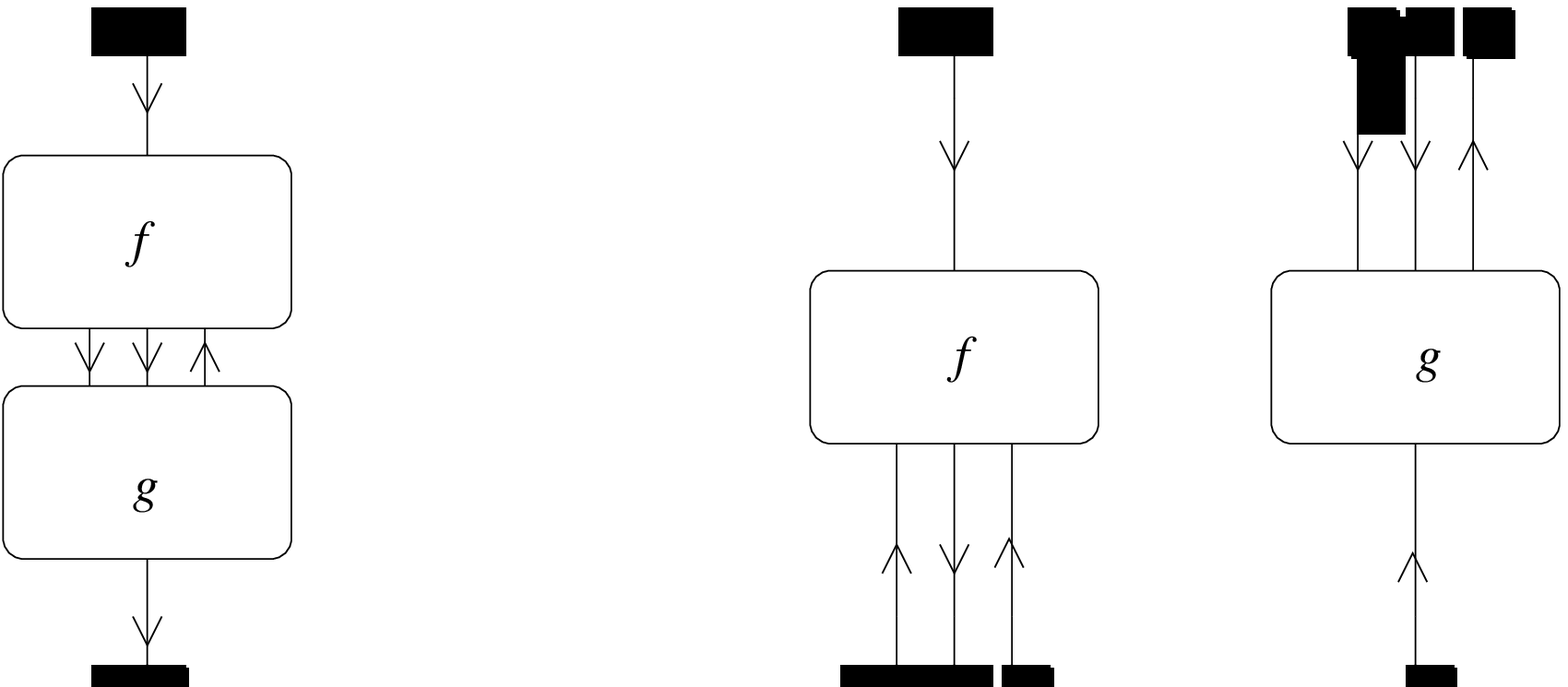}}
\medskip
\centerline{4.  Composition and tensor product of tangles}
\medskip
}
\noindent 
Thus it suffices to specify the morphisms in $H$ corresponding to certain
basic tangles from which all others can be built up by composition and
tensor product.  These basic tangles are shown in Figures 5 and 6.  A
downwards-pointing segment corresponds to the identity on $x$, while
an upwards-pointing segment corresponds to the identity on $x^\ast$.

\bigskip
\centerline{\epsfysize=0.75in\epsfbox{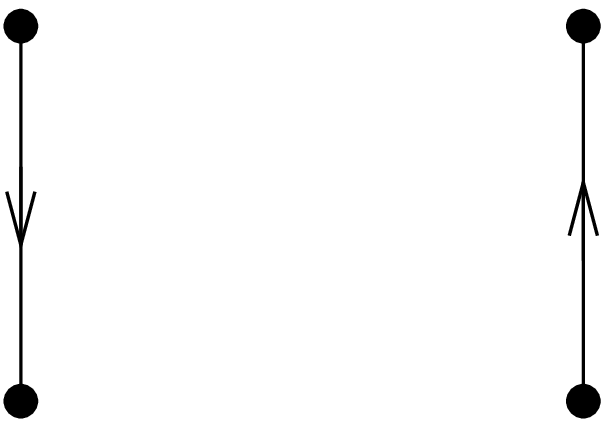}}
\medskip
\centerline{5.  Tangles corresponding to $1_x$ and $1_{x^\ast}$}
\medskip

\noindent
The two oriented forms of a `cup' tangle correspond to the morphisms
$e$ and $i^\ast$, while the two oriented forms of a `cap' correspond
to $i$ and $e^\ast$.

\bigskip
\centerline{\epsfysize=0.75in\epsfbox{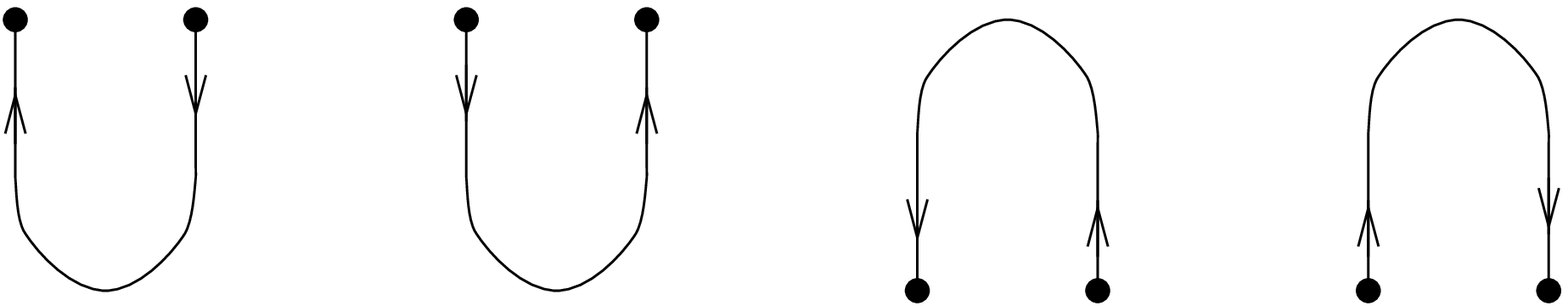}}
\medskip
\centerline{6. Tangles corresponding to $e$, $i^\ast$, $i$, and $e^\ast$}  
\medskip

\noindent
It then turns out that isotopic tangles correspond to the same
morphism in $H$.  The main thing to check is that the isotopic tangles
shown in Figure 7 correspond to the same morphisms.  This
follows from the triangle diagrams in the definition of an
adjunction.  Similar equations with the orientation of the arrows
reversed follow from Proposition \ref{2H*2}.

\bigskip
\centerline{\epsfysize=1.5in\epsfbox{triangle.eps}}
\medskip
\centerline{7. Tangle equations corresponding to the definition of adjunction}
\medskip

If $H$ is braided, we can also map framed oriented tangles in 3
dimensions to morphisms in $H$.  A typical such tangle is shown in
Figure 8.  We use the blackboard framing, in which each strand is
implicitly equipped with a vector field normal to the plane in which
the tangle is drawn.

\bigskip
\centerline{\epsfysize=1.5in\epsfbox{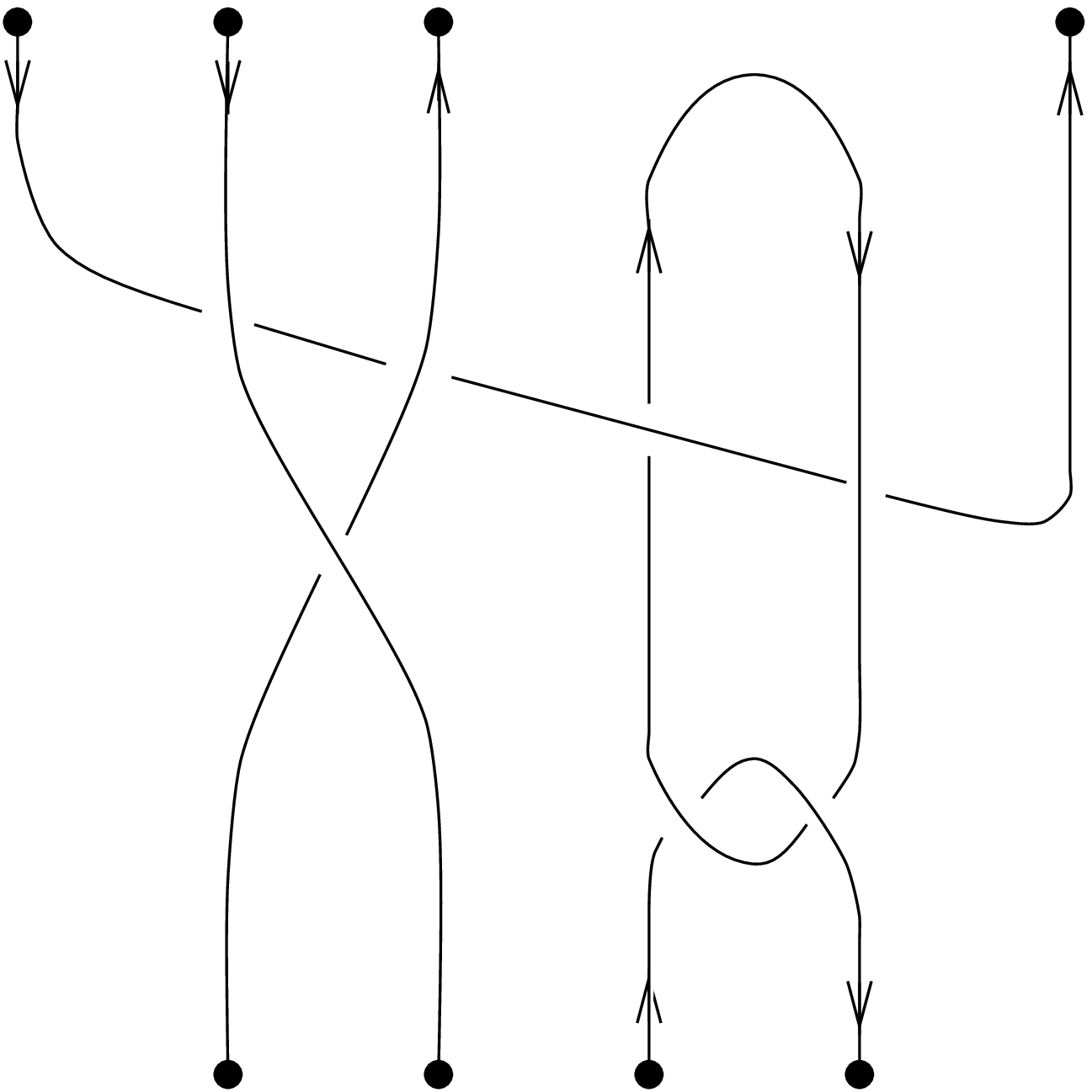}}
\medskip
\centerline{8. Typical tangle in 3 dimensions}
\medskip

\noindent  We interpret the basic tangles in Figures 5 and 6 as
we did before. Moreover, we let the tangles in Figure 9
correspond to the morphisms $B_{x,x}$, $B_{x^\ast,x}$,
$B_{x,x^\ast}$, and $B_{x^\ast,x^\ast}$, and let the tangles in
Figure 10 correspond to the morphisms $B_{x,x}^{-1}$,
$B_{x^\ast,x}^{-1}$, $B_{x,x^\ast}^{-1}$, and
$B_{x^\ast,x^\ast}^{-1}$.

\bigskip
\centerline{\epsfysize=1in\epsfbox{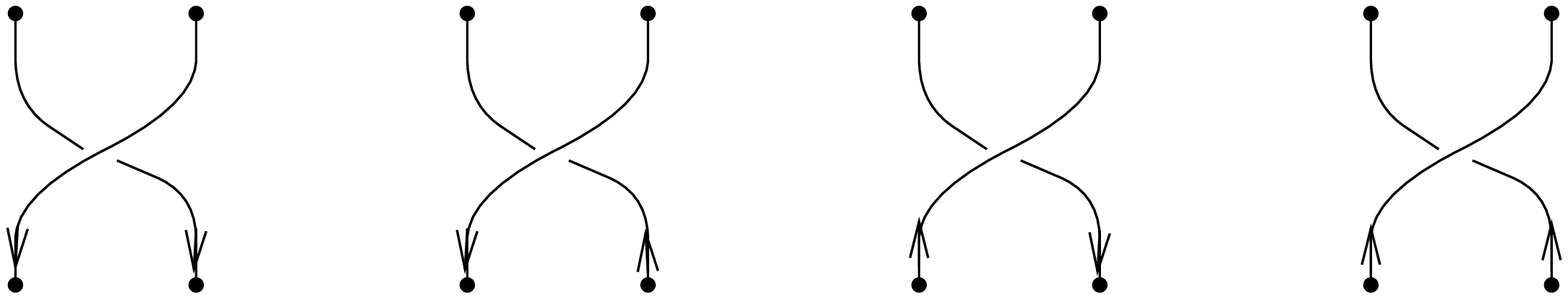}} 
\medskip
\centerline{7. Tangles corresponding to $B_{x,x}$, $B_{x^\ast,x}$,
$B_{x,x^\ast}$, and $B_{x^\ast,x^\ast}$}
\medskip
\bigskip
\centerline{\epsfysize=1in\epsfbox{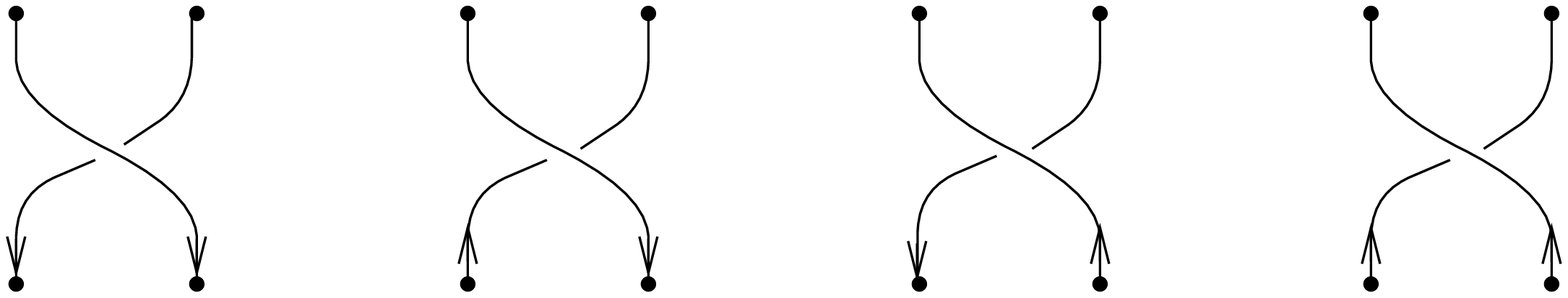}}
\medskip
\centerline{8. Tangles corresponding to $B_{x,x}^{-1}$, $B_{x^\ast,x}^{-1}$,
$B_{x,x^\ast}^{-1}$, and $B_{x^\ast,x^\ast}^{-1}$}
\medskip

Now suppose we wish isotopic framed oriented tangles to correspond to
the same morphism in $H$.  Invariance under the 2nd and 3rd
Reidemeister moves follows from the properties of the braiding,
so it suffices to check invariance under the framed version of the 1st
Reidemeister move.  For this, note that the tangle shown in
Figure 9 corresponds to the balancing of the adjunction $(x,x^\ast,i,e)$.  
This tangle has a $2\pi$ twist in its framing.

\vbox{
\bigskip
\centerline{\epsfysize=1.5in\epsfbox{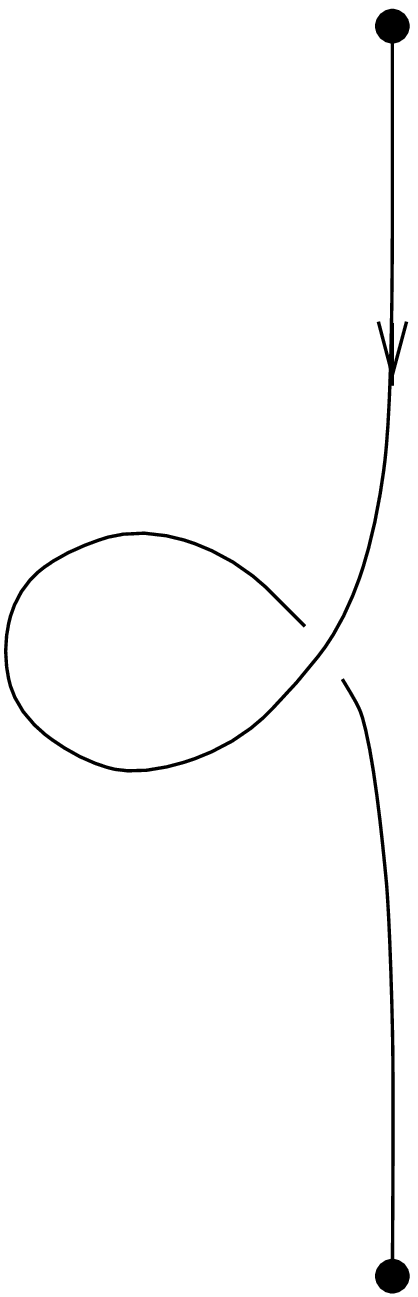}}
\medskip
\centerline{9. Tangle corresponding to the balancing $b \maps x \to x$}
\medskip
}
\noindent 
The framed version of the 1st Reidemeister move, shown in Figure 10,
represents the cancellation of two opposite $2\pi$ twists in the framing.
Both tangles in this picture correspond to the same morphism in $H$
precisely when the balancing $b \maps x \to x$ is unitary.  

\bigskip
\centerline{\epsfysize=1.5in\epsfbox{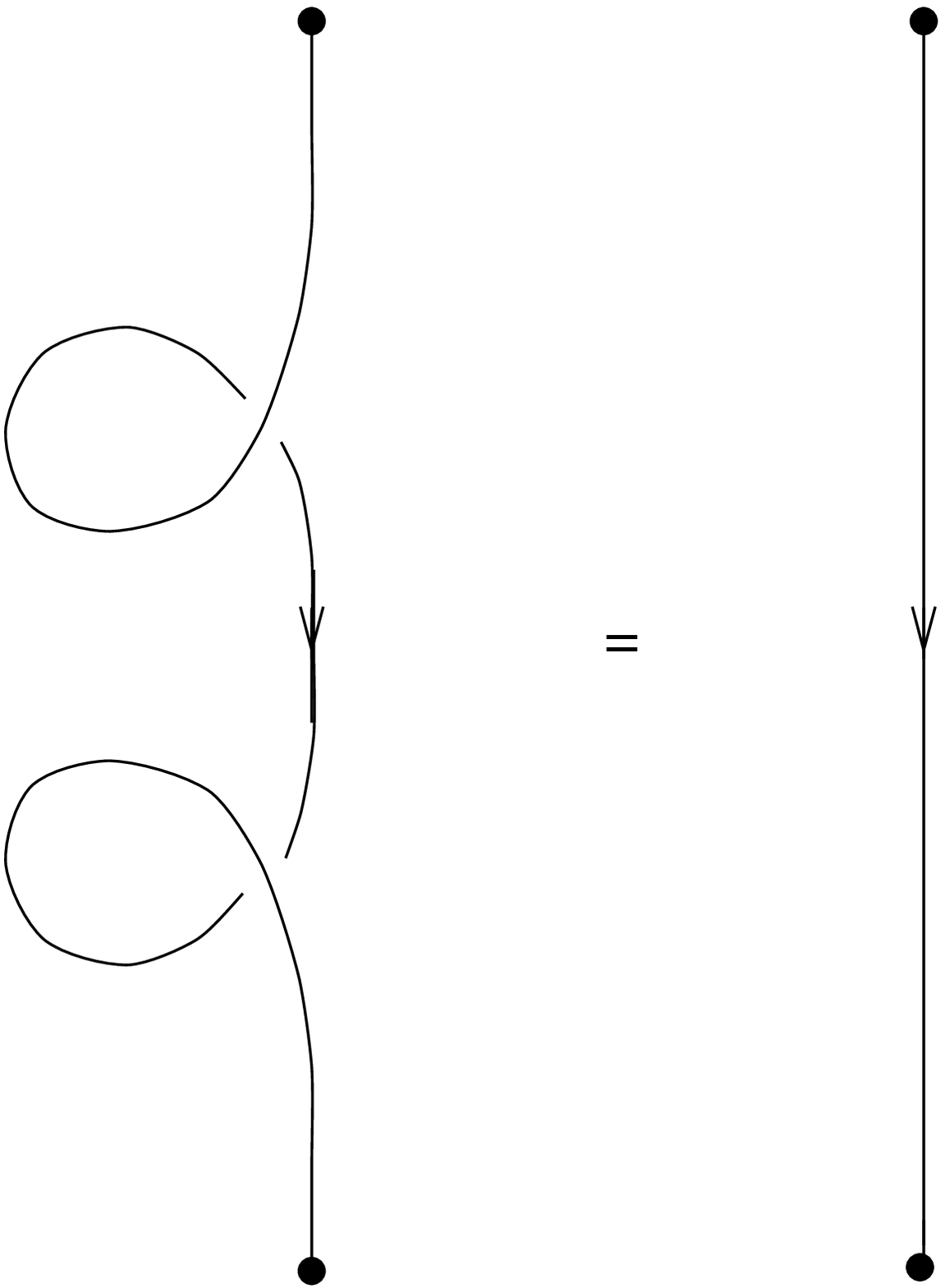}}
\medskip
\centerline{10. Tangle equation corresponding to unitarity of the
balancing}
\medskip

In short, we obtain a map from isotopy classes of framed oriented
tangles in 3 dimensions to morphisms in a braided 2-H*-algebra $H$
whenever we choose an adjunction in $H$ whose balancing is unitary.
This motivates the following definition:

\begin{defn} \et An adjunction $(x,x^\ast,i,e)$ in a braided 2-H*-algebra
is {\rm well-balanced} if its balancing is unitary.
\end{defn}

\noindent 
Similarly, given any well-balanced adjunction in a symmetric
2-H*-algebra $H$, we obtain a map from isotopy classes of framed
oriented tangles in 4 dimensions to morphisms in $H$.  We may draw
tangles in 4 dimensions just as we draw tangles in 3 dimensions, but
there is an extra rule saying that any right-handed crossing is
isotopic to the corresponding left-handed crossing.  One case of this
rule is shown in Figure 11.  Invariance under these isotopies follows
directly from the fact that the braiding is a symmetry.

\bigskip
\centerline{\epsfysize=1.0in\epsfbox{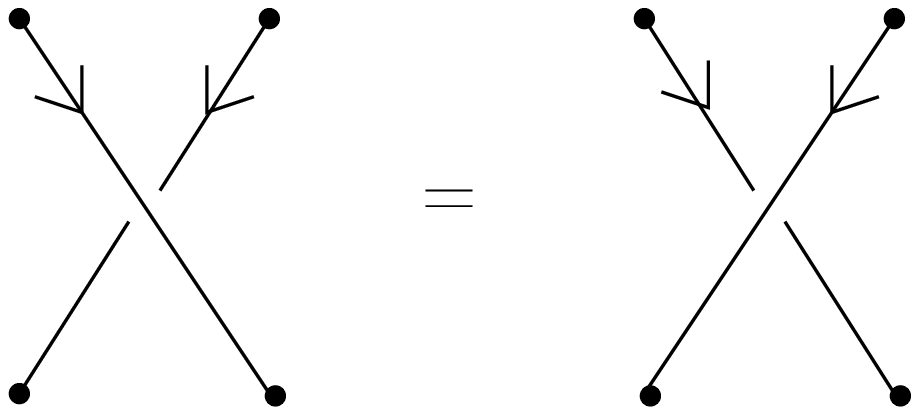}}
\medskip
\centerline{11. Tangle equation corresponding to symmetry}
\medskip

The important fact is that well-balanced adjunctions exist and
are unique up to a unique unitary isomorphism.  Moreover, all of
them have the same balancing:

\begin{thm}\label{2H*3} \et Suppose $H$ is a braided
2-H*-algebra.  For every object $x \in H$ there exists a well-balanced
adjunction $(x,y,i,e)$.  Given well-balanced adjunctions $(x,y,i,e)$
and $(x,y',i',e')$, there is a unique morphism $u \maps y \to y'$ such that 
\[ i' = i(u \tensor x), \qquad e' = (x \tensor u^{-1})e, \]
and this morphism is unitary.  \end{thm}

Proof - To simplify notation we assume without loss of generality that
$H$ is strict.  Suppose first that $x \in H$ is simple.  
Then for any adjunction $(x,y,i,e)$, the balancing equals $\beta 1_x$
for some nonzero $\beta \in \C$.  By Proposition \ref{2H*1} we may
define a new adjunction $(x,y,|\beta|^{1/2}i,|\beta|^{-1/2}e)$.
Since the balancing of this adjunction equals $\beta |\beta|^{-1} 1_x$, 
this adjunction is well-balanced.

Next suppose that $x \in H$ is arbitrary.  Using Theorem
\ref{H*3} we can write $x$ as an orthogonal direct sum of simple
objects $x_j$, in the sense that there are morphisms
\[            p_j \maps x \to x_j \]
with 
\[         p_j^\ast p_j = 1_{x_j}, \qquad \sum_j p_j p_j^\ast =
1_x.  \]
Let $y_j$ be a left dual of $x_j$, and
define $y$ to be an orthogonal direct sum 
of the objects $y_j$, with morphisms 
\[            q_j \maps y \to y_j \]
such that
\[         q_j^\ast q_j = 1_{y_j}, \qquad \sum_j q_j q_j^\ast =
1_y . \]
Since the $x_j$ are simple, there exist adjunctions
$(x_j,y_j,e_j,i_j)$ for which the balancings $b_j \maps x_j \to x_j$
are unitary.  Define the adjunction $(x,y,i,e)$ by
\[      i = \sum_j i_j (p_j^\ast \tensor q_j^\ast)  , \qquad
        e = \sum_j  (q_j \tensor p_j) e_j  .\]
One can check that this is indeed an adjunction and that 
the balancing $b \maps x \to x$ of this adjunction is given by
\[         b = \sum_j p_j b_j p_j^\ast ,\]
and is therefore unitary.  

Now suppose that $(x,y',i',e')$ is any other well-balanced adjunction with
$x$ as right dual.  Let $b'$ denote the balancing of this adjunction.
We shall prove that $b' = b$.
By Propositions \ref{H*4} and \ref{2H*1} there
exists a unitary morphism $g \maps y \to y'$, and we have
\ban  b' &=& (e'^\ast \tensor x)(y' \tensor B_{x,x})(e' \tensor x) \\
&=&  (e'^\ast(g^\ast \tensor x) \tensor x)(y \tensor B_{x,x})
((g \tensor x)e' \tensor x).  \ean
By Proposition \ref{2H*1}, $(x,y,(g \tensor x)e',i(x \tensor g^{-1}))$ is
an adjunction, so by the uniqueness up to isomorphism of right adjoints
we have $(g \tensor x)e' = (y \tensor f)e$ for some isomorphism
$f \maps x \to x$.  We thus have
\ban  b' &=&
(e^\ast(y \tensor f^\ast) \tensor x)(y \tensor B_{x,x})
((y \tensor f)e \tensor x) \\
&=& fbf^\ast .\ean
We may write $x$ as an orthogonal direct sum
\[        x = \bigoplus_\lambda x_\lambda \]
where $\{e_\lambda\}$ is a basis of $H$ and $x_\lambda$
is a direct sum of some number of copies of $e_\lambda$.
Then by our previous formula for $b$ we have
\[        b = \bigoplus_\lambda \beta_\lambda 1_{x_\lambda} \]
with $|\beta_\lambda| = 1$ for all $\lambda$.  
We also have
\[        f = \bigoplus_\lambda f_\lambda \]
for some morphisms $f_\lambda \maps x_\lambda \to x_\lambda$.  It follows that
\[        b' = \bigoplus_\lambda \beta_\lambda f_\lambda {f_\lambda}^\ast \]
Since $b$ and $b'$ are unitary it follows that each morphism
$f_\lambda {f_\lambda}^\ast$ is unitary.   Since the only positive
unitary operator is the identity, using Theorem \ref{H*3} it follows that 
each $f_\lambda {f_\lambda}^\ast$ is the identity,
so $b' = b$ as desired.  

By Proposition \ref{2H*1}, we know there is
a unique isomorphism $u \maps y \to y'$ with 
\[ i' = i(u \tensor x) , \qquad e' = (x \tensor u^{-1})e ,\]
and we need to show that $u$ is unitary.  Since $b' = b$, we have
\[ (ib' \tensor y)(x \tensor B_{y,y}^{-1})(i^\ast \tensor y)=
(ib \tensor y)(x \tensor B_{y,y}^{-1})(i^\ast \tensor y) ,\]
and if one simplifies this equation using the fact that
\[ b = (e^\ast \tensor x)(y' \tensor B_{x,x})(e \tensor x)  \]
and
\[ b' = (e^\ast(x \tensor (u^{-1})^\ast) \tensor x)(y' \tensor B_{x,x})
((x \tensor u^{-1})e \tensor x),\] 
one finds that $u$ is unitary.  \qed

\begin{cor} \et In a braided 2-H*-algebra every well-balanced
adjunction with  $x$ as right dual has the same balancing, which
we call {\rm the balancing} of $x$ and denote as $b_x \maps x \to
x$.  \end{cor}

Proof - This was shown in the proof above. \qed

Note that for any simple object $x$ in a braided 2-H*-algebra,
the balancing $b_x$ must equal $1_x$ times some unit complex number, the
{\it balancing phase} of $x$.  In physics, the balancing phase describes
the change in the wavefunction of a particle that undergoes a $2\pi$
rotation.   Note that in a symmetric 2-H*-algebra
\ban  b_x &=& (e_x^\ast \tensor 1_x)(1_{x^\ast} \tensor B_{x,x})
(e_x \tensor 1_x)  \\
&=&  (e_x^\ast \tensor 1_x)(1_{x^\ast} \tensor B^\ast_{x,x})
(e_x \tensor 1_x) \\
&=&  b_x^\ast, \ean
so $b_x^2 = 1_x$.  Thus in this case the balancing phase
of any simple object must be $\pm 1$.   In physics,
this corresponds to the fact that particles in 4-dimensional spacetime
are either bosons and fermions depending on the phase they acquire when
rotated by $2\pi$, while in 3-dimensional spacetime other 
possibilities, sometimes called `anyons', can occur \cite{DR,FK}.

More generally, we make the following definition:

\begin{defn} \et If $H$ is a symmetric 2-H*-algebra, an object $x \in H$ is
{\rm even} or {\rm bosonic} if $b_x = 1$, and {\rm odd} or {\rm
fermionic} if $b_x = -1$.  We say $H$ is {\rm even} or {\rm purely
bosonic} if every object $x \in H$ is even.  \end{defn} 

\noindent Note that if $x \oplus y$ is an orthogonal direct sum,
\[        b_{x \oplus y} = b_x \oplus b_y, \]
so an object in any symmetric 2-H*-algebra
is even (resp.\ odd) if and only if it is
a direct sum of even (resp.\ odd) simple objects.   Also, since 
\[       b_{x \tensor y} = (b_x \tensor b_y)B_{x,y}B_{y,x}, \]
it follows that the tensor product of two even or two odd objects is even, 
while the tensor product of an even and an odd object is odd.  

There is a way to turn any symmetric 2-H*-algebra into an even one,
which will be useful in Section \ref{recon}. 

\begin{prop} \label{bosonization}
\et {\rm (Doplicher-Roberts)}\, Suppose $H$ is a symmetric
2-H*-algebra.  Then there is a braiding $B'$ on $H$ given on simple
objects $x,y \in H$ by 
\[        B^\flat_{x,y} = (-1)^{|x|\,|y|} B_{x,y} \]
where $|x|$ equals $0$ or $1$ depending on whether $x$ is even or odd, and
similarly for $|y|$.   Let $H^\flat$ denote $H$ equipped with the
new braiding $B^\flat$.  Then $H^\flat$ is an even symmetric
2-H*-algebra, the {\rm bosonization} of $H$.  \end{prop} 

Proof - This is a series of straightforward computations.  One approach
involves noting that for any objects $x,y \in H$,
\[     B^\flat_{x,y} = {1\over 2}\, B_{x,y} (1_x \tensor 1_y + 1_x \tensor
b_y + b_x  \tensor 1_y - b_x \tensor b_y) . \]
\hbox{\hskip 30em} \qed

\noindent The above proposition is essentially due to Doplicher and
Roberts, who proved it in a slightly different context \cite{DR}.
However, the term `bosonization' is borrowed from Majid \cite{Majid},
who uses it to denote a related process that turns a super-Hopf algebra
into a Hopf algebra.

\subsection{Trace and dimension}

The notion of the `dimension' of an object in a braided
2-H*-algebra will be very important in Section \ref{recon}. 
First we introduce the related notion of `trace'.

\begin{defn}\et If $H$ is a braided 2-H*-algebra and $f \maps x \to x$ is a
morphism in $H$, for any well-balanced adjunction $(x,x^\ast,i,e)$
we define the {\rm trace} of $f$, $\tr(f) \in \End(1)$, by 
\[          \tr(f) = e (x^\ast \tensor f) e^\ast .\]
\end{defn}

\noindent
The trace is independent of the choice of well-balanced adjunction,
by Theorem \ref{2H*3}.  Also, one can show that an obvious
alternative definition of the trace is actually equivalent:
\[    \tr(f) = i^\ast (f \tensor x^\ast) i .\]

\begin{defn}\et If $H$ is a braided 2-H*-algebra, we define
the {\it dimension} of $x$, $\dim(x)$, to be $\tr(1_x)$.  \end{defn}

\noindent Note that $x,y$ are objects in a braided 2-H*-algebra, we have
\[       \dim(x \oplus y) = \dim(x) + \dim(y), \qquad
         \dim(x \tensor y) = \dim(x) \dim(y), \qquad
          \dim(x^\ast) = \dim(x) .\]
Moreoever, we have:

\begin{prop} \label{dim} \et If $H$ is a symmetric 2-H*-algebra and $x
\in H$ is any object, then the spectrum of $\dim(x)$ is a subset of $\N
= \{0,1,2,\dots\}$. 
\end{prop}

Proof - We follow the argument of Doplicher and Roberts \cite{DR}.  
For any $n \ge 0$, the group algebra of the symmetric group $S_n$ acts
as endomorphisms of $x^{\tensor n}$, and the morphisms $p_S, p_A$
corresponding to complete symmetrization and complete
antisymmetrization, respectively, are self-adjoint projections in the
H*-algebra $\End(x^{\tensor n})$.   It follows that 
$\tr(p_S), \tr(p_A) \ge 0$.  If $x$ is even, a calculation shows that
\[    \tr(p_A) = {1\over n!}\, \dim(x)(\dim(x) - 1)\cdots (\dim(x) - n +
1) \]
For this to be nonnegative for all $n$, the spectrum of $\dim(x)$ must
lie in $\N$.  Similarly, if $x$ is odd, a calculation shows that 
\[    \tr(p_S) = {1\over n!} \, \dim(x)(\dim(x) - 1)\cdots (\dim(x) - n +
1) \]
so again the spectrum of $\dim(x)$ lies in $\N$.  In general, any
object $\dim(x)$ is a sum of simple objects, which are either even or
odd, so by the additivity of dimension, the spectrum of $\dim(x)$ again
lies in $\N$.  \qed

For any 2-H*-algebra, the Eckmann-Hilton argument shows that $\End(1)$
is a commutative H*-algebra, and thus isomorphic to a direct sum of
copies of $\C$.  (See HDA0 or HDA1 for a explanation of the Eckmann-Hilton
argument.)  

\begin{defn} \et  A 2-H*-algebra $H$ is {\rm connected} if the unit
object $1 \in H$ is simple.  \end{defn}

\noindent In a connected 2-H*-algebra, $\End(1) \iso \C$.  The
dimension of any object in a connected symmetric 2-H*-algebra
is thus a nonnegative integer.  

In addition to the above notion of dimension it is also interesting
to consider the `quantum dimension'.  Here our treatment most closely
parallels that of Majid \cite{Majid}.  

\begin{defn} \et 
If $H$ is a braided 2-H*-algebra and $f \maps x \to x$ is a
morphism in $H$, for any well-balanced adjunction $(x,x^\ast,i,e)$
we define the {\rm quantum trace} of $f$, $\qtr(f) \in \End(1)$, by 
\[          \qtr(f) = \tr(b_x f) .\]
We define the {\rm quantum dimension} of $x$, $\qdim(x)$, to be
$\qtr(1_x)$.  \end{defn}

\noindent In the case of a symmetric 2-H*-algebra, the quantum trace
is also called the `supertrace'.  Suppose $H$ is a connected symmetric 
2-H*-algebra and $x$ is a simple object.  Then $\qdim(x) \ge 0$ if
$x$ is even and $\qdim(x) \le 0$ if $x$ is odd.  The idea of odd
objects as negative-dimensional is implicit in Penrose's work 
on negative-dimensional vector spaces \cite{Penrose}.  

\subsection{Homomorphisms and 2-homomorphisms}

There is a 2-category with 2-H*-algebras as objects and `homomorphisms' and
`2-homomorphisms' as morphisms and 2-morphisms, respectively.  This
is also true for braided 2-H*-algebras and symmetric 2-H*-algebras.

\begin{defn} \et Given 2-H*-algebras $H$ and $H'$, a {\rm
homomorphism} $F \maps H \to H'$ is a morphism of 2-Hilbert spaces
that is also a monoidal functor.  If $H$ and $H'$ are braided, we say
that $F$ is a {\rm homomorphism of braided 2-H*-algebras} if $F$ is
additionally a braided monoidal functor.  If $H$ and $H'$ are symmetric,
we say that $F$ is a {\rm homomorphism of symmetric 2-H*-algebras} if $F$ is
a morphism of 2-Hilbert spaces that is also a symmetric monoidal functor.
\end{defn}

Recall that a functor $F \maps C \to C'$ between monoidal categories is
monoidal if it is equipped with a natural isomorphism
$\Phi_{x,y} \maps F(x) \tensor F(y) \to F(x \tensor y)$
making the following diagram commute for any objects $x,y,z \in C$:
\[
\begin{diagram} [F(x) \tensor (F(y) \tensor F(z))]
\node{(F(x) \tensor F(y)) \tensor F(z)}
\arrow{e,t}{\Phi_{x,y}\, \tensor 1_{F(z)}}  \arrow{s,l}{a_{F(x),F(y),F(z)}} 
\node{F(x \tensor y) \tensor F(z)} 
\arrow{e,t}{\Phi_{x\tensor y,z}} 
\node{F((x \tensor y)\tensor z)} 
\arrow{s,r}{F(a_{x,y,z})}  \\
\node{F(x) \tensor (F(y) \tensor F(z))} 
\arrow{e,t}{1_{F(x)} \tensor \Phi_{y,z}}
\node{F(x) \tensor F(y \tensor z)}
\arrow{e,t}{\Phi_{x,y\tensor z}}
\node{F(x \tensor (y \tensor z))}
\end{diagram}
\]
together with an isomorphism $\phi \maps 1_{C'} \to
F(1_{C})$ making 
the following diagrams commute for any object $x \in C$:
\[
\begin{diagram}[F(1) \tensor F(x)]
\node{1 \tensor F(x)} \arrow{e,t}{\ell_{F(x)}} 
\arrow{s,l}{\phi \tensor 1_{F(x)}}
\node{F(x)} \\
\node{F(1) \tensor F(x)} \arrow{e,t}{\Phi_{1,x}} 
\node{F(1 \tensor x)} \arrow{n,r}{F(\ell_x)}
\end{diagram}
\]
\[  
\begin{diagram}[F(1) \tensor F(x)]
\node{F(x) \tensor 1} \arrow{e,t}{r_{F(x)}} 
\arrow{s,l}{1_{F(x)}\tensor \phi}
\node{F(x)} \\
\node{F(x) \tensor F(1)} \arrow{e,t}{\Phi_{x,1}} 
\node{F(x \tensor 1)} \arrow{n,r}{F(r_x)}
\end{diagram}
\]
If $C$ and $C'$ are braided, we say that $F$ is braided if additionally
it makes the following diagram commute for all $x,y \in C$:
\[  
\begin{diagram}[F(x) \tensor F(y)]
\node{F(x) \tensor F(y)} \arrow{e,t}{B_{F(x),F(y)}} 
\arrow{s,l}{\Phi_{x,y}}
\node{F(y) \tensor F(x)}
\arrow{s,r}{\Phi_{y,x}} \\
\node{F(x \tensor y)} \arrow{e,t}{F(B_{x,y})} 
\node{F(y \tensor x)} 
\end{diagram}
\]
A symmetric monoidal functor is simply a braided monoidal functor that
happens to go between symmetric monoidal categories!  No extra condition
is involved here.  

Note that if $F \maps H \to H'$ is a homomorphism of braided 2-H*-algebras,
$F$ maps any well-balanced adjunction in $H$ to one in $H'$.  Thus it
preserves dimension in the following sense:
\[           \dim(F(x)) = F(\dim(x))  \]
for any object $x \in H$.
In particular, if $H$ and $H'$ are connected, so that we can
identify the dimension of objects in either with numbers, we have
simply $\dim(F(x)) = \dim(x)$.  

\begin{defn}\et If $H$ and $H'$ are 2-H*-algebras, possibly braided or
symmetric, and $F,G \maps H \to H'$ are homomorphisms of the appropriate
sort, a {\rm 2-homomorphism} $\alpha \maps F \To G$ is a monoidal
natural transformation.  \end{defn}

Suppose that the $(F,\Phi,\phi)$ and $(G,\Gamma,\gamma)$ are monoidal
functors from the monoidal category $C$ to the monoidal category $D$. 
Then a natural transformation $\alpha \maps F \to G$ is monoidal if 
the diagrams
\[  \begin{diagram}[F(x) \tensor F(y)]
\node{F(x) \tensor F(y)} \arrow{e,t} {\alpha_x \tensor \alpha_y}
\arrow{s,l}{\Phi_{x,y}}  
\node{G(x) \tensor G(y)} \arrow{s,r}{\Gamma_{x,y}}  
\\ \node{F(x \tensor y)}  \arrow{e,t}{\alpha_{x\tensor y}} 
\node{G(x \tensor y)}  
\end{diagram}
\]
and
\[      \begin{diagram}[F(1)] 
\node{1} \arrow{s,l}{\phi} \arrow{se,t}{\gamma} \\
\node{F(1)} \arrow{e,t}{\alpha_1} \node{G(1)} 
\end{diagram}  
\]
commute.  
There are no extra conditions required of `braided monoidal' or
`symmetric monoidal' natural transformations.  

Finally, when we speak of two 2-H*-algebras $H$ and $H'$, possibly
braided or symmetric, being {\it equivalent}, we always mean the
existence of homomorphisms $F \maps H \to H'$ and $G \maps H' \to H$ of the
appropriate sort that are inverses up to a 2-isomorphism.

\section{Reconstruction Theorems} \label{recon}

In this section we give a classification of symmetric 2-H*-algebras.
Doplicher and Roberts proved a theorem which implies that connected
even symmetric 2-H*-algebras are all equivalent to categories of compact
group representations \cite{DR,DR2}.  Here and in all that follows, by a
`representation' of a compact group we mean a finite-dimensional
continuous unitary representation.  Given a compact group $G$, let
$\Rep(G)$ denote the category of such representations of $G$.  This
becomes a connected even symmetric 2-H*-algebra in an obvious way.
While Doplicher and Roberts worked using the language of `C*-categories',
their result can be stated as follows:

\begin{thm}\label{dhr1}\et  {\rm (Doplicher-Roberts)} Let $H$ be
a connected even symmetric 2-H*-\break algebra.  Then there
exists a homomorphism of symmetric 2-H*-algebras $T \maps H \to
\Hilb$, unique up to a unitary 2-homomorphism.  Let $U(T)$ be the
group of unitary 2-homomorphisms $\alpha \maps T \To T$, given
the topology in which a net $\alpha_\lambda \in U(T)$ converges
to $\alpha$ if and only if $(\alpha_\lambda)_x \to \alpha_x$ in
norm for all $x \in H$.  Then $U(T)$ is compact, each Hilbert
space $T(x)$ becomes a representation of $U(T)$, and the
resulting homomorphism $\tilde T \maps H \to \Rep(U(T))$ extends
to an equivalence of symmetric 2-H*-algebras.  \end{thm}

Note that any continuous homomorphism $\rho \maps G \to G'$ between
compact groups determines a homomorphism of symmetric 2-H*-algebras,
\[        \rho^\ast \maps \Rep(G') \to \Rep(G), \]
sending each representation $\sigma$ of $G'$ to the representation
$\sigma \circ \rho$ of $G$.  The above theorem yields a useful converse
to this construction:

\begin{cor}\label{dhr2}\et {\rm (Doplicher-Roberts)} Let $F \maps H' \to
H$ be a homomorphism of connected even symmetric 2-H*-algebras.  Let
$T \maps H \to \Hilb$ be a homomorphism of
symmetric 2-H*-algebras.  Then there exists a continuous group
homomorphism  
\[            F^\ast \maps U(T) \to U(FT) \]
such that $F^\ast(\alpha)$ equals the horizontal composite $F\circ \alpha$.
Moreover, $(F^\ast)^\ast$ equals $F$ up to a unitary 2-homomorphism.
\end{cor} 

Dolan \cite{Dolan} has noted that a generalization of the
Doplicher-Roberts theorem to even symmetric 2-H*-algebras --- not
necessarily connected --- amounts to a categorification of the
Gelfand-Naimark theorem.  The spectrum of a commutative H*-algebra $H$
is a set $\Spec(H)$ whose points are homomorphisms from $H$ to $\C$.
The Gelfand-Naimark theorem implies that $H$ is isomorphic to the
algebra of functions from $\Spec(H)$ to $\C$.  Similarly, we may
define the `spectrum' of an even symmetric 2-H*-algebra $H$ to be the
groupoid $\Spec(H)$ whose objects are homomorphisms from $H$ to
$\Hilb$, and whose morphisms are unitary 2-homomorphisms between
these.  Moreover, we shall show that $H$ is equivalent to a symmetric
2-H*-algebra whose objects are `representations' of $\Spec(H)$ ---
certain functors from $\Spec(H)$ to $\Hilb$.  Indeed, our proof of
this uses an equivalence
\[      \hat{\hbox{\hskip 0.5em}}\maps H \to \Rep(\Spec(H)) \]
that is just the categorified version of the `Gelfand transform' 
for commutative H*-algebras.  

In fact, there is no need to restrict ourselves to symmetric
2-H*-algebras that are even.  To treat a general symmetric
2-H*-algebra $H$ we need objects of $\Spec(H)$ to be
homomorphisms from $H$ to a symmetric 2-H*-algebra of
`super-Hilbert spaces'.   The spectrum will then be a
`supergroupoid' --- though not the most general sort of thing one
could imagine calling a supergroupoid.

\begin{defn} \et Define $\SuperHilb$ to be the category whose
objects are $Z_2$-graded (finite-dimensional) Hilbert spaces, 
and whose morphisms are linear maps preserving the grading.  \end{defn}

The category $\SuperHilb$ can be made into a symmetric 2-H*-algebra
where the $\ast$-structure is the ordinary Hilbert space adjoint, 
the product is the usual tensor product of $\Z_2$-graded
Hilbert spaces, and the braiding is given on homogeneous elements $v \in x$, 
$w \in y$ by
\[           B_{x,y}(v \tensor w) = (-1)^{{\rm deg} v\, {\rm \deg} w}
\, w \tensor v.  \]

\begin{defn} \et If $H$ is a symmetric 2-H*-algebra, define $\Spec(H)$
to be the category whose objects are symmetric 2-H*-algebra homomorphisms
$F \maps H \to \SuperHilb$ and whose morphisms are unitary
2-homomorphisms between these.  \end{defn}

\begin{defn} \et A {\rm topological groupoid} is a groupoid for
which the $\hom$-sets are topological spaces and
the groupoid operations are continuous.  A {\rm compact groupoid} is
is a topological groupoid with compact Hausdorff $\hom$-sets and finitely
many isomorphism classes of objects. \end{defn}

\begin{defn} \et A {\rm supergroupoid} is a groupoid $G$ equipped
with a natural transformation $\beta \maps 1_G \To 1_G$, the {\rm
balancing,} with $\beta^2 = 1$. A {\rm compact supergroupoid} is a
supergroupoid that is also a compact groupoid.  \end{defn}

Let $H$ be a symmetric 2-H*-algebra.  Then $\Spec(H)$ becomes a
topological groupoid if for any $S,T \maps H \to \Hilb$ we give
$\hom(S,T)$ the topology in which a net $\alpha_\lambda$ converges to
$\alpha$ if and only if $(\alpha_\lambda)_x \to \alpha_x$ in norm for
any $x \in H$.  We shall show that $\Spec(H)$ is a compact groupoid.
Also, $\Spec(H)$ becomes a supergroupoid if for any object $T \in
\Spec(H)$ we define $\beta_T \maps T \To T$ by
\[        (\beta_T)_x = b_{T(x)} = T(b_x)  \]
for any object $x \in H$.  One can
check that $\beta \maps 1_{\Spec(H)} \To 1_{\Spec(H)}$ is a natural
transformation, and $\beta^2 = 1$ because the balancing for $H$ satisfies
$b_x^2 = 1$ for any $x \in H$.  

\begin{defn}\et Given a compact supergroupoid $G$, a (continuous, unitary,
finite-dimensional) {\rm representation} of $G$ is a functor $F \maps G
\to \SuperHilb$ such that $F(g)$ is unitary for every morphism $g$ in $G$, 
$F \maps \hom(x,y) \to \hom(F(x),F(y))$ is continuous for all objects
$x,y \in G$, and 
\[         F(\beta_x) = b_{F(x)}  \]
for every object $x \in G$.  We define $\Rep(G)$ to be the category having
representations of $G$ as objects and natural transformations between
these as morphisms.  \end{defn} 

Let $G$ be a compact supergroupoid.  Then the category $\Rep(G)$
becomes an even symmetric 2-H*-algebra in a more or less obvious way
as follows.  Given objects $F,F' \in \Rep(G)$, we make $\hom(F,F')$
into a Hilbert space with the obvious linear structure and the inner
product given by
\[        \langle \alpha, \beta \rangle = \sum_x 
\tr({\alpha_x}^\ast \beta_x) \]
where the sum is taken over any maximal set of nonisomorphic objects of
$G$.   This makes $\Rep(G)$ into a $\Hilb$-category.
Moreover, $\Rep(G)$ becomes a 2-Hilbert space if we define the dual
of $\alpha \maps F \To F'$ by $(\alpha^\ast)_x = (\alpha_x)^\ast$. 
We define the tensor product of objects $F,F' \in \Rep(G)$ by
\[       (F \tensor F')(x) = F(x) \tensor F'(x), \qquad
         (F \tensor F')(f) = F(f) \tensor F'(f) \] 
for any object $x \in G$ and morphism $f$ in $G$.  It is easy to define a
tensor product of morphisms and associator making $\Rep(G)$ into a
monoidal category, and to check that $\Rep(G)$ is then a 2-H*-algebra.
Finally, $\Rep(G)$ inherits a braiding from the braiding in $\SuperHilb$,
making $\Rep(G)$ into a symmetric 2-H*-algebra.

Now suppose $H$ is an even symmetric 2-H*-algebra.  Then there is a functor
\[      \hat{\hbox{\hskip 0.5em}}\maps H \to \Rep(\Spec(H)) , \]
the {\it categorified Gelfand transform}, given as follows.
For every object $x \in H$, $\hat x$ is the representation with
\[        \hat x(T) = T(x) \]
for all $T \in \Spec(H)$, and
\[        \hat x(\alpha) = \alpha_x \]
for all $\alpha \maps T \To T'$, where $T,T' \in \Spec(H)$.  For every
morphism $f \maps x \to y$ in $H$, $\hat f \maps \hat x \To \hat y$ is
the natural transformation with 
\[    \hat f(T) = T(f) \]
for all $T \in \Spec(H)$.  Our generalized Doplicher-Roberts theorem
states: 

\begin{thm} \label{dhr4}\et Suppose that $H$ is a
symmetric 2-H*-algebra.  Then $\Spec(H)$ is a compact supergroupoid
and $\hat{\hbox{\hskip 0.5em}}\maps H \to \Rep(\Spec(H))$ extends to
an equivalence of symmetric 2-H*-algebras.
\end{thm}

Proof - We have described how $\Spec(H)$ is a supergroupoid.    To see
that it is compact, note that for any $S,T \in  \Spec(H)$ the $\hom$-set
$\hom(S,T)$ is a compact Hausdorff space, by Tychonoff's theorem.  We
also need to show that $\Spec(H)$ has finitely many isomorphism classes
of objects.   The unit object $1_H$ is the direct sum of finitely many
nonisomorphic simple objects $e_i$, the kernels of the minimal
projections $p_i$ in the commutative H*-algebra $\End(1_H)$.   Any
object $x \in H$ is thus a direct sum of objects $x_i = e_i \tensor x$,
and any morphism $f \maps x \to y$ is a direct sum of morphisms $f_i
\maps x_i \to y_i$.  In short, $H$ is, in a fairly obvious sense, the
direct sum of finitely many connected symmetric 2-H*-algebras $H_i$.  
Any homomorphism $T \maps H \to \SuperHilb$ induces a homomorphism from
$\End(1_H)$ to $\End(1_{\SuperHilb}) \iso \C$, which must annihilate all
but one of the projections $p_i$, so $T$ sends one of the objects $x_i$
to $1_{\SuperHilb}$ and the rest to $0$.  Thus $\Spec(H)$ is, as a
groupoid, equivalent to the disjoint union of the groupoids
$\Spec(H_i)$, and hence has finitely many isomorphism classes of objects.

To show that the categorified Gelfand transform is an equivalence, first
suppose that $H$ is even and connected.  Then the supergroupoid
$\Spec(H)$ has $\beta = 1$, so every representation $F \maps \Spec(H)
\to \SuperHilb$ factors through the inclusion $\Hilb \hookrightarrow
\SuperHilb$.  Moreover, by Theorem \ref{dhr1} all the objects of
$\Spec(H)$ are isomorphic, so $\Spec(H)$ is equivalent, as a groupoid, 
to the group $\U(T)$ for any $T \in \Spec(H)$.  We thus
obtain an equivalence of symmetric 2-H*-algebras between
$\Rep(\Spec(H))$ and $\Rep(\U(T))$ as defined in Theorem \ref{dhr1}. 
Using this, the fact that $\tilde T \maps H \to \Rep(\U(T))$ is
an equivalence translates into the fact that $\hat{\hbox{\hskip
0.5em}}\maps H \to \Rep(\Spec(H))$ is an equivalence.

Next, suppose that $H$ is even but not connected.  Then $H$ is a direct
sum of the even connected symmetric 2-H*-algebras $H_i$ as above, and
$\Rep(\Spec(H))$ is similarly the direct sum of the $\Rep(\Spec(H_i))$.  
Because the categorified Gelfand transform $\hat{\hbox{\hskip
0.5em}}\maps H_i \to \Rep(\Spec(H_i))$ is an equivalence for all $i$, 
$\hat{\hbox{\hskip 0.5em}}\maps H \to \Rep(\Spec(H))$ is an equivalence. 

Finally we treat the general case where $H$ is an arbitrary
symmetric 2-H*-algebra.  Note that if $H$ and $K$ are
symmetric 2-H*-algebras, a symmetric 2-H*-algebra homomorphism $F
\maps H \to K$ gives rise to a symmetric 2-H*-algebra
homomorphism $F^\flat \maps H^\flat \to K^\flat$ between their
bosonizations, where $F^\flat$ is the same as $F$ on objects and
morphisms.  Note also that $F$ is an equivalence of symmetric
2-H*-algebras if and only if $F^\flat$ is.  Thus to show that
$\hat{\hbox{\hskip 0.5em}}\maps H \to \Rep(\Spec(H))$ is an
equivalence, it suffices to show $\hat{\hbox{\hskip 0.5em}}^\flat
\maps H^\flat \to \Rep(\Spec(H))^\flat$ is an equivalence.

For this, note that any supergroupoid $G$ has a {\it
bosonization} $G^\flat$, in which the underlying compact groupoid
of $G$ is equipped with the trivial balancing $\beta = 1$.  Moreover,
there is a homomorphism of symmetric 2-H*-algebras
\[   
\begin{diagram}[\Rep(G)^\flat]
\node{\Rep(G)^\flat} \arrow{e,t}{X} \node{\Rep(G^\flat)} 
\end{diagram} 
\]
sending any representation $F \in \Rep(G)^\flat$ to the representation 
$X(F) \in \Rep(G^\flat)$ given by the commutative square
\[  
\begin{diagram}[\SuperHilb]
\node{G^\flat} \arrow{e,t}{X(F)} \arrow{s,l}{I} \node{\SuperHilb} \\
\node{G} \arrow{e,t}{F} \node{\SuperHilb} \arrow{n,r}{E}
\end{diagram}
\]
Here $I \maps G^\flat \to G$ is the identity on the underlying
groupoids, while the 2-H*-algebra homomorphism $E \maps \SuperHilb \to
\SuperHilb$ maps any super-Hilbert space to the even super-Hilbert space
with the same underlying Hilbert space, and acts as the identity on
morphisms.  One may check that $X(F)$ is really a compact supergroupoid
representation.  Similarly, given a morphism $\alpha \maps F \To F'$ in
$\Rep(G)^\flat$, we define $X(\alpha)$ to be the horizontal composite $I
\circ \alpha \circ E$. In fact, $X$ is an equivalence, for given any
representation $F$ of $G^\flat$ we can turn it back into a
representation of $G$ by equipping each Hilbert space $F(x)$, $x \in G$
with the grading $F(\beta_x)$, where $\beta$ is the balancing of $G$.

Similarly, for any symmetric 2-H*-algebra $H$ there is an equivalence
\[
\begin{diagram}[\Spec(H)^\flat]
\node{\Spec(H)^\flat} \arrow{e,t}{Y} \node{\Spec(H^\flat)} 
\end{diagram} 
\]
sending any object $T \in \Spec(H)^\flat$ to the object 
$Y(T) \in \Spec(H^\flat)$ given by the commutative square
\[  
\begin{diagram}[\SuperHilb]
\node{H^\flat} \arrow{e,t}{Y(T)} \arrow{s,l}{I} \node{\SuperHilb} \\
\node{H} \arrow{e,t}{T} \node{\SuperHilb} \arrow{n,r}{E}
\end{diagram}
\]
where $I \maps H^\flat \to H$ is the identity on the underlying
2-H*-algebras, while $E$ is given as above.

We thus have equivalences
\[ \begin{diagram}[\Rep(\Spec(H))^\flat]
\node{\Rep(\Spec(H))^\flat} \arrow{e,t}{\sim}
\node{\Rep(\Spec(H)^\flat)} \arrow{e,t}{\sim}
\node{\Rep(\Spec(H^\flat))} 
\end{diagram}
\]
and their composite gives a diagram commuting up to natural isomorphism:
\[ \begin{diagram}[\Rep(\Spec(H))^\flat]
\node{H^\flat} \arrow{e,t}{\hat{\hbox{\hskip 0.5em}}^\flat}
\arrow{se,b}{\hat{\hbox{\hskip 0.5em}}}
\node{\Rep(\Spec(H))^\flat} \arrow{s}\\
\node[2]{\Rep(\Spec(H^\flat))}  
\end{diagram}
\]
It follows that $\hat{\hbox{\hskip 0.5em}}^\flat
\maps H^\flat \to \Rep(\Spec(H))^\flat$ is an equivalence, as was to 
be shown.  \hbox{\hskip 30em} \qed

Presumably what Theorem \ref{dhr4} is trying to tell us is
that there are 2-functors $\Rep$ and $\Spec$ going both ways
between the 2-category of compact supergroupoids and the 2-category
of symmetric 2-H*-algebras, and that these extend to a 2-equivalence
of 2-categories.  We shall not try to prove this here.  However, it
is worth noting that for any compact supergroupoid $G$, there is
a functor 
\[  \check{\hbox{\hskip 0.5em}} \maps G \to \Spec(\Rep(G)) \]
given as follows.  For every object $x \in G$, $\hat x$ is the
object of $\Spec(\Rep(G))$ with
\[             \hat x(F) = F(x)  \]
for all $F \in \Rep(G)$, and
\[              \hat x(\alpha) = \alpha_x \]
for all $\alpha \maps F \to F'$, where $F,F' \in \Rep(G)$.  For
every morphism $g \maps x \to y$ in $G$, $\check g \maps \check x \To
\check y$ is the natural transformation with 
\[       \check g(F) = F(g) \]
for all $F \in \Rep(G)$.  Presumably
$\check{\hbox{\hskip 0.5em}} \maps G \to \Spec(\Rep(G))$ is in some
sense an equivalence of compact supergroupoids.

\subsection{Compact abelian groups}

The representation theory of compact abelian groups is rendered
especially simple by the use of Fourier analysis, as generalized by
Pontryagin.  Suppose that $T$ is a compact abelian group.  Then its
dual $\hat T$ is defined as the set of equivalence classes of
irreducible representations $\rho$ of $T$.  The dual becomes a
discrete abelian group with operations given as follows:
\[           [\rho][\rho'] = [\rho \tensor \rho'], \]
\[          [\rho]^{-1} = [\rho^\ast]  \]
Then the Fourier transform is a unitary isomorphism 
\[            f \maps L^2(T) \to L^2(\hat T) \]
given by 
\[            f(\chi_\rho) = \delta_{[\rho]} \]
where $\chi_\rho$ is the character of the representation $\rho$,
and $\delta_{[\rho]}$ is the function on $\hat T$ which equals
$1$ at $[\rho]$ and $0$ elsewhere.

The Fourier transform has an interesting categorification. Note that
the ordinary Fourier transform has as its domain the
infinite-dimensional Hilbert space $L^2(T)$, which has a basis given
by the characters of irreducible representations of $T$.  The
categorified Fourier transform will have as domain the 2-Hilbert space
$\Rep(T)$, which has a basis given by the irreducible representations
themselves.  (Taking the character of a representation is a form of
`decategorification'.)  Similarly, just as the ordinary Fourier
transform has as its codomain an infinite-dimensional Hilbert space of
$\C$-valued functions on $\hat T$, the categorified Fourier transform
will have as its codomain a 2-Hilbert space of $\Hilb$-valued
functions on $\hat T$.

More precisely, define $\Hilb[G]$ for any discrete group $G$ to be the
category whose objects are $G$-graded Hilbert spaces for which the
total dimension is finite, and whose morphisms are linear maps
preserving the grading.  Alternatively, we can think of $\Hilb[G]$ as the
category of hermitian vector bundles over $G$ for which the sum of the
dimensions of the fibers is finite.  We may write any object $x \in
\Hilb[G]$ as a $G$-tuple $\{x(g)\}_{g \in G}$ of Hilbert spaces.  The
category $\Hilb[G]$ becomes a 2-H*-algebra in an obvious way with a
product modelled after the convolution product in the group algebra $\C[G]$:
\[          (x \tensor y)(g) = \bigoplus_{\{g',g'' \in G\,\colon\; g'g'' = g\}}
x(g') \tensor y(g'')  .\]
If $G$ is abelian, $\Hilb[G]$ becomes a symmetric 2-H*-algebra.

Now suppose that $T$ is a compact abelian group.  Given any object $x
\in \Rep(T)$, we may decompose $x$ into subspaces corresponding to the
irreducible representations of $T$:
\[        x = \bigoplus_{g \in \hat T} x(g) . \]
We define the {\it categorified Fourier transform} 
\[        F \maps \Rep(T) \to \Hilb[\hat T]  \]
as follows.  For any object $x \in \Rep(T)$, we set
\[        F(x) = \{x(g)\}_{g \in \hat T} . \]
Moreoever, any morphism $f \maps x \to y$ in $\Rep(T)$ gives rise to
linear maps $f(g) \maps x(g) \to y(g)$ and thus a morphism
$F(f)$ in $\Hilb[\hat T]$.  One can check that $F$ is not only 2-Hilbert
space morphism but actually a homomorphism of symmetric 2-H*-algebras.
This is the categorified analog of how the ordinary Fourier transform
sends pointwise multiplication to convolution.  Note that, in analogy
to the formula 
\[            f(\chi_\rho) = \delta_{[\rho]} \]
satisfied by the ordinary Fourier transform, for any irreducible
representation $\rho$ of $T$ the categorified Fourier transform
$F(\rho)$ is a hermitian vector bundle that is 1-dimensional at
$[\rho]$ and 0-dimensional elsewhere.

\begin{thm} \et If $T$ is a compact abelian group, the
categorified Fourier transform $F \maps \Rep(T) \to \Hilb(\hat
T)$ is an equivalence of symmetric 2-H*-algebras.  \end{thm}

Proof - There is a homomorphism $G \maps \Hilb[\hat T] \to
\Rep(T)$ sending each object $\{x(g)\}_{g \in \hat T}$ in
$\Hilb[\hat T]$ to a representation of $T$ which is a direct sum
of spaces $x(g)$ transforming according to the different
isomorphism classes $g \in \hat T$ of irreducible representations
of $T$.  One can check that $FG$ and $GF$ are naturally
isomorphic to the identity.  \qed

%\subsection{Compact Lie groups}

%\begin{defn} \label{generated} \et We say a 2-H*-algebra $H$ is
%{\rm generated} by a set $S$ of objects if any 2-H*-subalgebra
%containing the objects in $S$ is equivalent, as a symmetric
%2-H*-algebra, to $H$.  \end{defn}

%\begin{thm} \et Let $H$ be a connected even symmetric 2-H*-algebra.
%Then the following are equivalent:

%\begin{enumerate}
%\item $H$ is generated by a single object.  
%\item $H$ is generated by a finite set of objects.  
%\item $H$ is equivalent as a symmetric 2-H*-algebra to $\Rep(G)$ for
%some compact Lie group $G$.  
%\end{enumerate}
%\end{thm}

\subsection{Compact classical groups} \label{classical}

The representation theory of a `classical' compact Lie group has a
different flavor from that of general compact Lie groups.  The
representation theory of general compact Lie groups heavily involves
the notions of maximal torus, Weyl group, roots and weights.  We hope
to interpret this theory in terms of 2-Hilbert spaces in a future
paper.  However, the representation theory of a classical group can
also be studied using Young diagrams \cite{Weyl}.  This approach
relies on the fact that its categories of representations have simple
universal properties.  These universal properties can be described in
the language of symmetric 2-H*-algebras, and a description along these
lines represents a distilled version of the Young diagram theory.

For example, consider the group $\U(n)$.  The fundamental
representation of $\U(n)$ on $\C^n$ is the `universal $n$-dimensional
representation'.  In other words, for a group to have a
(unitary) representation on $\C^n$ is precisely for it to have a
homomorphism to $\U(n)$.  This universal property can also be
expressed as a universal property of $\Rep(\U(n))$.  Suppose that $G$
is a compact group.  Then any $n$-dimensional representation $y \in
\Rep(G)$ is isomorphic to a representation of the form $\rho \maps G
\to \U(n)$.  The representation $\rho$ gives rise to a homomorphism
\[       \rho^\ast \maps \Rep(\U(n)) \to \Rep(G),   \]
and letting $x$ denote the fundamental representation of $\U(n)$, we have
$\rho^\ast(x) = \rho$.  Since $\rho$ and $y$ are isomorphic, there is a
unitary 2-homomorphism from $\rho^\ast$ to a homomorphism
\[       F \maps \Rep(\U(n)) \to \Rep(G) \]
with $F(x) = y$.  

In short, for any $n$-dimensional object $y \in \Rep(G)$ there is a
homomorphism $F \maps \Rep(\U(n)) \to \Rep(G)$ of symmetric
2-H*-algebras with $F(x) = y$.  On the other hand, suppose $F' \maps
\Rep(\U(n)) \to \Rep(G)$ is any other homomorphism with $F'(x) = y$.  We
claim that there is a unitary 2-homomorphism from $F$ to $F'$.  By
Corollary \ref{dhr2}, there exists a homomorphism $\rho' \maps G \to
\U(n)$ with a unitary 2-homomorphism from $F'$ to $\rho'^\ast$.  On the
other hand, by construction there is a unitary 2-homomorphism from $F$
to $\rho^\ast$ for some $\rho \maps G \to \U(n)$.  To show there is a
unitary 2-homomorphism from $F$ to $F'$, it thus suffices to show that
$\rho$ and $\rho'$ are isomorphic in $\Rep(G)$.  This holds because
$\rho \iso y = F'(x) \iso \rho'^\ast(x) = \rho'$.

Now, since any connected even symmetric 2-H*-algebra is unitarily
equivalent to $\Rep(G)$ for some compact $G$ by Theorem \ref{dhr1}, we
may restate these results as follows.  Suppose $H$ is a connected even
symmetric 2-H*-algebra and let $y$ be an $n$-dimensional object of $H$. 
Then there exists a homomorphism $F \maps \Rep(\U(n)) \to H$ with $F(x)
= y$.  Moreover, this is unique up to a unitary 2-homomorphism. 
Furthermore, we can drop the assumption that $H$ is even by working with
the full subcategory whose objects are all the even objects of $H$.

We may thus state the universal property of $\Rep(\U(n))$ as
follows: 

\begin{thm}\label{un}\et $\Rep(\U(n))$ is the free connected symmetric
2-H*-algebra on an even object $x$ of dimension $n$.  That is, given any
even $n$-dimensional object $y$ of a connected symmetric 2-H*-algebra
$H$, there exists a homomorphism of symmetric 2-H*-algebras $F \maps
\Rep(U(n)) \to H$ with $F(x) = y$, and $F$ is unique up to a unitary
2-homomorphism.  
\end{thm}  

Let $\Lambda^n x$ denote the cokernel of $p_A \maps x^{\tensor n} \to
x^{\tensor n}$ (complete antisymmetrization), and let $S^n x$ denote
the cokernel of $p_S \maps x^{\tensor n} \to x^{\tensor n}$ (complete
symmetrization).  We can describe the category of representations of
$\SU(n)$ as follows:

\begin{thm}\label{sun} \et $\Rep(\SU(n))$ is the free connected symmetric 
2-H*-algebra on an even object $x$ with $\Lambda^n x \iso 1$.  \end{thm}

Proof - Suppose that $G$ is a compact group and the object $y \in
\Rep(G)$ has $\Lambda^n y \iso 1$.  It follows that $y$ is
$n$-dimensional by the computation in Proposition \ref{dim}, and the
isomorphism $\Lambda^n y \iso 1$ determines a $G$-invariant volume
form on the representation $y$.  Thus $y$ is isomorphic to a
representation of the form $\rho \maps G \to \SU(n)$.  The rest of the
proof follows that of Theorem \ref{un}. \qed

Here we can see in a simple context how our theory is a distillation
of the theory of Young diagrams.  (The Young diagram approach to 
representation theory is more familiar for $\SU(n)$ than for $\U(n)$.)
In heuristic terms, the above theorem says that every representation of
$\SU(n)$ is generated from the fundamental representation $x$ using
the operations present in a symmetric 2-H*-algebra --- the
$\ast$-structure, direct sums, cokernels, tensor products, duals, and
the symmetry --- with no relations other than those implied by the
axioms for a connected symmetric 2-H*-algebra and the fact that
$x$ is even and $\Lambda^n x \iso 1$.  The theory of Young diagrams
makes this explicit by listing the irreducible representations of
$\SU(n)$ in terms of minimal projections $p \maps x^{\tensor x} \to
x^{\tensor k}$, or in other words, Young diagrams with $k$ boxes.  The
symmetric 2-H*-algebra of representations of a subgroup $G \subset
\SU(n)$, such as $\SO(n)$ or $\Sp(n)$, is a quotient of
$\Rep(\SU(n))$.  We may describe this quotienting process by giving
extra relations as in Theorems \ref{spn} and \ref{son} below.  These
extra relations give identities saying that different Young diagrams
correspond to the same representation of $G$.

The classical groups $\O(n)$ and $\Sp(n)$ are related to the concept
of self-duality.  Given adjunctions $(x,x^\ast,i_x,e_x)$ and
$(y,y^\ast,i_y,e_y)$ in a monoidal category $C$, for any morphism $f
\maps x \to y$ there is a morphism $f^\dagger \maps y^\ast \to
x^\ast$, given in the strict case by the composite:
\[ \begin{diagram}[y^\ast \tensor x \tensor x^\ast]   
\node{y^\ast = y^\ast \tensor 1} \arrow{e,t}{y^\ast \tensor i_x}
\node{y^\ast \tensor x \tensor x^\ast} 
\arrow{e,t}{y^\ast \tensor f \tensor x^\ast}
\node{y^\ast \tensor y \tensor x^\ast} \arrow{e,t}{e_y \tensor x^\ast}
\node{1 \tensor x^\ast = x^\ast}
\end{diagram}
\]  
(Our notation here differs from that of HDA0.)  Since the left dual of
an object in a 2-Hilbert space is also its right dual as in
Proposition \ref{2H*2}, given a morphism $f \maps x \to x^\ast$ we
obtain another morphism $f^\dagger \maps x \to x^\ast$.  Using this we
may describe $\Rep(\O(n))$ and $\Rep(\Sp(n))$ as certain `free connected
symmetric 2-H*-algebras on one self-dual object':

\begin{thm}\label{on} \et  $\Rep(\O(n))$ is the free connected symmetric
2-H*-algebra on an even object $x$ of dimension $n$ with an
isomorphism $f \maps x \to x^\ast$ such that $f^\dagger =
f$. \end{thm}

Proof - 
Suppose that $G$ is a compact group and the object $y \in
\Rep(G)$ is $n$-dimensional and equipped with an isomorphism
$f \maps x \to x^\ast$ with $f^\dagger = f$.  Then there is a
nondegenerate pairing $F \maps y \tensor y \to 1$ given by $F = (y
\tensor f)i_y^\ast$.  A calculation, given in the proof of Proposition
\ref{selfdual}, shows that $F$ is symmetric.  It follows that $y$ is
isomorphic to a representation of the form $\rho \maps G \to \O(n)$.
The rest of the proof follows that of Theorem \ref{un}. \qed

\begin{thm}\label{spn} \et $\Rep(\Sp(n))$ is the free connected symmetric
2-H*-algebra on one even object $x$ with $\Lambda^n x \iso 1$ 
and with an isomorphism $f \maps x \to x^\ast$ such that $f^\dagger =
-f$.  \end{thm}

Proof - The proof is analogous to that of Theorem \ref{on}, except
that the pairing $F$ is skew-symmetric.  \qed

Following the proof of Theorem \ref{sun} we may also characterize
$\Rep(\SO(n))$ as follows:

\begin{thm}\label{son} \et $\Rep(\SO(n))$ is the free connected symmetric
2-H*-algebra on an even object $x$ with $\Lambda^n x \iso 1$  
and with an isomorphism $f \maps x \to x^\ast$ such that $f^\dagger =
f$.  \end{thm} 

The conditions on the isomorphism $f \maps x \to x^\ast$ in Theorems
\ref{on} and \ref{spn} are quite reasonable, in the following sense:

\begin{prop}\label{selfdual} \et 
Suppose that $x$ is a simple object in a symmetric 
2-H*-algebra and that $x$ is isomorphic to $x^\ast$.  Then one and only
one of the following is true: either there is an isomorphism
$f \maps x \to x^\ast$ with $f = f^\dagger$, or there is an isomorphism
$f \maps x \to x^\ast$ with $f = -f^\dagger$.  \end{prop}

Proof - Note that there is an isomorphism of complex vector spaces
\ban    \hom(x,x^\ast) &\iso& \hom(x \tensor x,1) \\
                  f &\mapsto& (1 \tensor f)i_x^\ast, \ean
and note that 
\[         \hom(x \tensor x,1) \iso
\hom(S^2 x,1) \oplus \hom(\Lambda^2 x,1) . \]
Suppose $f \maps x \to x^\ast$ is an
isomorphism and let $F = (1 \tensor f)i_x^\ast$.  Since $x$ is simple,
$f$ and thus $F$ is unique up to a scalar multiple, so $F$
must lie either in $\hom(S^2 x,1)$ or $\hom(\Lambda^2 x,1)$. 
In other words, $B_{x,x}F = \pm F$.  Choose a well-balanced
adjunction for $x$.  Assuming without loss of generality that $H$ is 
strict, we have 
\ban  f^\dagger 
&=& (x \tensor i_x)(x \tensor f\tensor x^\ast)(i_x^\ast \tensor x^\ast) \\
&=& (x \tensor i_x)(F \tensor x^\ast) \\
&=&\pm (x \tensor i_x)(B_{x,x}F \tensor x^\ast) \\
&=&\pm (f \tensor i_x)(B_{x^\ast,x} \tensor x^\ast)(i_x^\ast \tensor x^\ast)\\
&=& \pm f b_{x^\ast} \ean
Since $b_{x^\ast} = \pm 1_{x^\ast}$ depending on whether $x$,
and thus $x^\ast$, is even or odd, we have $f^\dagger = \pm f$.  \qed

\noindent This result is well-known if $H$ is a category of compact
group representations \cite{FH}.  Here one may also think of the morphism
$f \maps x \to x^\ast$ as a conjugate-linear intertwining operator
$j \maps x \to x$.  The condition that $f = \pm f^\dagger$ is then
equivalent to the condition that $j^2 = \pm 1_x$.  
One says that $x$ is a real representation if $j^2 = 1_x$ and 
a quaternionic representation if $j^2 = -1_x$, establishing the
useful correspondence:
\[   {\rm real:complex:quaternionic::orthogonal:unitary:symplectic} \]

The following alternate characterization of $\Rep(\U(1))$ is
interesting because it emphasizes the relation between duals and
inverses.  Whenever $T$ is a compact abelian group and $x \in
\Rep(T)$, the dual $x^\ast$ is also the inverse of $x$, in the sense
that $x \tensor x^\ast \iso 1$.  We have:

\begin{thm}\et $\Rep(\U(1))$ is the free connected symmetric 2-H*-algebra
on an even object $x$ with $x \tensor x^\ast \iso 1$.  \end{thm} 

Proof - By Theorem \ref{un} it suffices to show that an object $x$ in a
connected symmetric 2-H*-algebra is 1-dimensional if and only if $x
\tensor x^\ast \iso 1$.  On the one hand, by the multiplicativity of
dimension, $x \tensor x^\ast \iso 1$ implies that $\dim(x) = 1$.  On the
other hand, suppose $\dim(x) = 1$.  Then we claim $i_x \maps 1 \to x
\tensor x^\ast$ and $i_x^\ast \maps x \tensor x^\ast \to 1$ are inverses.
First, ${i_x}^\ast i_x$ is the identity since $\dim(x) = 1$.  Second,
${i_x}^\ast i_x \in \End(x \tensor x^\ast)$ is idempotent since $\dim(x) =
1$.  Since $x \tensor x^\ast$ is 1-dimensional, it is simple (by the
additivity of dimension), so ${i_x}^\ast i_x$ must be the identity. \qed

Finally, it is interesting to note that $\SuperHilb$ is the free connected 
symmetric 2-H*-algebra on an odd object $x$ with $x \tensor x \iso 1$.
This object $x$ is the one-dimensional odd super-Hilbert space.

\section{Conclusions}

The reader will have noted that some of our results are slight reworkings of
those in the literature.  One advantage of our approach is that it
immediately suggests generalizations to arbitrary $n$.  While the
general study of $n$-Hilbert spaces will require a deeper
understanding of $n$-category theory, we expect many of the same
themes to be of interest.  With this in mind, let us point out some
problems with what we have done so far.

One problem concerns the definition of the quantum-theoretic hierarchy.
A monoid is a essentially a category with one object.  More precisely, a
category with one object $x$ can be reconstructed from the monoid
$\End(x)$, and up to isomorphism every monoid comes from a one-object
category in this way.  Comparing Figures 1 and 2, one might at first
hope that by analogy an H*-algebra would be a one-dimensional 2-Hilbert
space.  Unfortunately, the way we have set things up, this is not the case.

If $H$ is a one-dimensional 2-Hilbert space with basis given by
the object $x$, then $\End(x)$ is an H*-algebra.   However,
$\End(x)$ is always isomorphic to $\C$; one does not get any
other H*-algebras this way.  The reason appears to be the
requirement that a 2-Hilbert space has cokernels, so that if
$\End(x)$ has nontrivial idempotents, $x$ has subobjects.  If we
dropped this clause in the definition of a 2-Hilbert space, there
would be a correspondence between H*-algebras and 2-Hilbert
spaces, all of whose objects are direct sums of a single object
$x$. Perhaps in the long run it will be worthwhile to modify the
definition of 2-Hilbert space in this way.   On the other hand,
an H*-algebra is also an H*-category with one object.  An
H*-category has sums and differences of morphisms, but not of
objects, i.e., it need not have direct sums and cokernels. 
Perhaps, therefore, a $k$-tuply monoidal $n$-Hilbert space should
really be some sort of `$(n+k)$-H*-category' with one
$j$-morphism for $j < k$, and sums and differences of
$j$-morphisms for $j \ge k$.

A second problem concerns the program of getting an invariant of
$n$-tangles in $(n+k)$-dimensions from an object in a $k$-tuply monoidal
$n$-Hilbert space.  Let us recall what is known so far here.

Oriented tangles in 2 dimensions are the morphisms in a monoidal
category with duals, $C_{1,1}$.  Here by `monoidal category with duals'
we mean a monoidal $\ast$-category in which every object has a left
dual, the tensor product is a $\ast$-functor, and the associator is a
unitary natural transformation.  Suppose that $X$ is any other monoidal
category with duals, e.g.\ a 2-H*-algebra.  Then any adjunction
$(x,x^\ast,i,e)$ in $C$ uniquely determines a monoidal $\ast$-functor $F
\maps C_{1,1} \to H$ up to monoidal unitary natural isomorphism.  The
functor $F$ is determined by the requirement that it maps the positively
oriented point to $x$, the negatively oriented point to $x^\ast$, and
the appropriately oriented `cup' and `cap' tangles to $e$ and $i$.  

According to our philosophy we would prefer $F$ to be determined by an
object $x \in X$ rather than an adjunction.  However $F$ is not
determined up to natural transformation by requiring that it map the
positively oriented point to $x$.  For example, take $X = \Hilb$ and let $x
\in X$ be any object.   We may let $F$ send the negatively oriented point
to the dual Hilbert space $x^\ast$, and send the cap and cup to the
standard linear maps $e \maps x^\ast \tensor x \to \C$ and $i \maps \C
\to x \tensor x^\ast$.  Then $F(ii^\ast) = \dim(x)1_x$.  Alternatively
we may let $F$ send the cap and cup to $e' = \alpha^{-1} e$ and $i' =
\alpha i$ for any nonzero $\alpha \in \C$.  Then $F(ii^\ast) =
|\alpha|^2 \dim(x)1_x$.   The problem is that while adjunctions in $X$
are unique up to unique isomorphism, the isomorphism is not
necessarily unitary.  

In HDA0 we outlined a way to deal with this problem by
`strictifying'  the notion of a monoidal category with duals. 
Roughly speaking, this amounts to equipping each object with a
choice of left adjunction, and requiring the functor $F \maps
C_{1,1} \to X$ to preserve this choice.   Then $F$ is determined
up to monoidal unitary natural transformation by the requirement
that it map the positively oriented point to a particular object
$x \in X$.  In this paper we have attempted to take the `weak'
rather than the `strict' approach.  Our point here is that the
weak approach seems to make it more difficult to formulate the
sense in which $C_{1,1}$ is the `free' monoidal category with
duals on one object.

In higher dimensions the balancing plays an interesting role in this
issue.  Framed oriented tangles in 3 dimensions form a braided monoidal
category $C_{1,2}$ with duals.  Here by `braided monoidal category with
duals' we mean a monoidal category with duals which is also braided,
such that the braiding is unitary and every object $x$ has a
well-balanced adjunction $(x,x^\ast,i,e)$.    For any object $x$ in a
braided monoidal category $X$ with duals,  there is a braided monoidal
$\ast$-functor $F \maps C_{1,2} \to X$ sending the positively oriented
point to $x$.  Moreover, because well-balanced adjunctions are unique
up to unique {\it unitary} isomorphism, $F$ is unique up to monoidal
unitary natural isomorphism.  This gives a sense in which $C_{1,2}$
is the free braided monoidal category with duals on one object.  

Similarly, framed oriented tangles in 4 dimensions form a symmetric
monoidal category with duals $C_{1,3}$, i.e., a braided monoidal
category with duals for which the braiding is a symmetry.  Again, for
any object $x$ in a symmetric monoidal category $X$ with duals, there is
a symmetric monoidal $\ast$-functor $F \maps C_{1,3} \to X$ sending the
positively oriented point to $x$, and $F$ is unique up to monoidal
unitary natural isomorphism.  (For an alternative `strict' approach to
the 3- and 4-dimensional cases, see HDA0.)

In short, we need to understand the notion of $k$-tuply monoidal
$n$-Hilbert spaces more deeply, as well as the notion of `free'
$k$-tuply monoidal $n$-categories with duals.  

\subsection*{Acknowledgements}
Many of the basic ideas behind this paper were developed in collaboration
with James Dolan.  I would also like to thank Louis Crane, Martin Hyland,
Martin Neuchl, John Power, Stephen Sawin, Gavin Wraith, and David Yetter for
helpful conversations and correspondence.  I am grateful to 
the Erwin Schr\"odinger Institute and the physics department of
Imperial College for their hospitality while part of this work was
being done.

\end{document}